%% file: alternating-trees.tex
\PassOptionsToPackage{dvipsnames}{xcolor}

\documentclass{amsart}

\usepackage[T1]{fontenc} 
\usepackage[utf8]{inputenc}

\usepackage[colorlinks,citecolor=blue,urlcolor=blue]{hyperref}

\usepackage[backend=bibtex,doi=true,isbn=false,url=false,style=alphabetic]{biblatex}
\usepackage{amsthm,amsmath,amsfonts,amssymb}
\usepackage{commath} 
\usepackage{nicefrac} 
\usepackage{mathdots}
\usepackage{bbm}

\usepackage[capitalize,noabbrev]{cleveref} 

\usepackage[backgroundcolor=yellow]{todonotes} 

\usepackage{subfiles} 

\usepackage{float}
\usepackage{xspace} 
\usepackage{enumitem} 
\usepackage[config,labelfont={normalsize}]{caption,subfig}

\usepackage{xcolor}

\usepackage{tikz} 
	\usetikzlibrary{arrows.meta} 
	\usetikzlibrary{colorbrewer}
	\usetikzlibrary{patterns}
	\usetikzlibrary{matrix,arrows,decorations.pathmorphing}
	\usetikzlibrary{shapes.geometric} 
	\usetikzlibrary{decorations.markings}
	\usetikzlibrary{calc}
	\tikzstyle{very densely dashed}=[dash pattern=on 5pt off 1.5pt] 
	\tikzstyle{very very densely dashed}=[dash pattern=on 8pt off 1.5pt] 
	\tikzstyle{middle densely dashed}=[dash pattern=on 4pt off 1.5pt] 
	\tikzstyle{space densely dashed}=[dash pattern=on 0pt off 0.4pt on 10pt off 0pt]

\usepackage{graphicx}

\numberwithin{equation}{section}

\theoremstyle{plain} 
\newtheorem{thm}{Theorem}[section]
\newtheorem{theorem}[thm]{Theorem} 

\newtheorem{lemma}[thm]{Lemma} 
\newtheorem{problem}[thm]{Problem}
\newtheorem{corollary}[thm]{Corollary}
\newtheorem{proposition}[thm]{Proposition}

\theoremstyle{remark} 
\newtheorem{remark}[thm]{Remark}
\newtheorem{example}[thm]{Example}

\Crefname{equation}{}{Equations}

\makeatletter \newcommand{\@giventhatstar}[2]{\left[#1\;\middle|\;#2\right]}
\newcommand{\@giventhatnostar}[3][]{#1(#2\;#1|\;#3#1)}
\newcommand{\giventhat}{\@ifstar\@giventhatstar\@giventhatnostar} \makeatother

\newcommand{\macierz}{A} 
\newcommand{\g}{H} 
\newcommand{\cauchy}{\mathbf{G}}
\newcommand{\sS}{S} 
\newcommand{\gS}{F} 
\DeclareMathOperator{\sMS}{MS}
\newcommand{\gMS}{F} %
\newcommand{\sC}{C} 
\newcommand{\gC}{G} 
\newcommand{\gbi}{G} %
\DeclareMathOperator{\sMC}{MC} 
\newcommand{\gMC}{G}
\newcommand{\gMCempty}{G_{\emptyset}} %
\newcommand{\sT}{\mathbb{T}} 
\newcommand{\gT}{H} 
\newcommand{\sV}{V}
\newcommand{\sB}{B} 
\newcommand{\sW}{W} 
\newcommand{\sR}{R} 
\newcommand{\sE}{E}

\newcommand{\E}{\mathbb{E}} 
\newcommand{\R}{\mathbb{R}}
\newcommand{\N}{\mathbb{N}} 
\newcommand{\Z}{\mathbb{Z}}
 
\newcommand{\Moment}{M}
\newcommand{\Pro}{\mathbb{P}} 

\newcommand{\indicator}{\mathbbm{1}}

\newcommand{\nklatki}{n} 
\newcommand{\nmoment}{k} 
\newcommand{\h}{a}
\newcommand{\en}{\ell} 
\newcommand{\te}{\ell} 

\newcommand{\grafydek}{\mathcal{G}}

 \newcommand{\Sy}[1]{\mathfrak{S}_{#1}}
\newcommand{\func}{\Theta}

\newcommand{\weight}{w} 
\DeclareMathOperator{\comp}{Comp}
\DeclareMathOperator{\com}{Comp}

\DeclareMathOperator{\both}{tm} 
\DeclareMathOperator{\tagg}{tag} 
\DeclareMathOperator{\markk}{mark}
\DeclareMathOperator{\labeled}{lab} 

\DeclareMathOperator{\Ins}{Ins}
\DeclareMathOperator{\uIns}{u-Ins} 
\DeclareMathOperator{\Var}{Var}

\DeclareMathOperator{\RSK}{RSK}

\newcommand{\F}{\mathcal{F}}

\newcommand{\Deco}{\mathbbm{X}}

\newcommand{\kerx}{\mathbbm{x}} 
\newcommand{\kery}{\mathbbm{y}}
\newcommand{\kerell}{\mathbb{L}}

\newcommand{\sniady}{the second named author\xspace}

\newcommand{\iid}{i.i.d.\@\xspace}

\newcommand{\simplex}{\mathcal{S}}

\newcommand{\diagram}{\xi} 

\newcommand{\yd}[1]{\Xi_{#1}}

\input{figures/preamble-figures.tex}

\begin{document}

\title[Cumulants of threshold for Schensted insertion]{Cumulants of threshold \\ for Schensted row
    insertion \\ into random tableaux}

\author{Mikołaj Marciniak} 
\address{Interdisciplinary Doctoral School ``Academia Copernicana'', Faculty of
Mathematics and Computer Science, \mbox{Nicolaus Copernicus University in
    Toru{\'n}}, ul.~Fryderyka Chopina 12/18, 87-100~Toru{\'n}, Poland}
\email{marciniak@int.pl}

\author{Piotr Śniady} 
\address{Institute of
    Mathematics, Polish Academy of Sciences, \mbox{ul.~\'Sniadec\-kich 8,} 00-656~Warszawa, Poland}
\email{psniady@impan.pl} 

\begin{abstract}
    \emph{Schensted row insertion} is a fundamental component of
    the \emph{Ro\-bin\-son--Schensted--Knuth (RSK) algorithm}, a powerful tool in
    combinatorics and representation theory. This study examines the insertion of a
    deterministic number into a random tableau of a specified shape, focusing on the
    relationship between the value of the inserted number and the position of the
    new box created by the Schensted row insertion. Specifically, for a given
    tableau and a point on its boundary, we consider the threshold that separates
    values which, if inserted, would result in the new box being created above the
    point from those that would result in a new box below. We analyze a random
    tableau of fixed shape and study the corresponding random threshold value.
    Explicit combinatorial formulas for the cumulants of this random variable are
    provided, expressed in terms of \emph{Kerov's transition measure of the
        diagram}. These combinatorial formulas involve summing over \emph{non-crossing
        alternating trees}. As a first application of these results, we demonstrate that
    for random Young tableaux of prescribed large shape, the rightmost entry in the
    first row converges in distribution to an explicit Gaussian distribution.
\end{abstract}



\maketitle

\subfile{SECTION-teaser.tex}

\subfile{SECTION-general-form.tex}

\subfile{SECTION-cumulative-function}

\subfile{SECTION-trees.tex}

\section{Application:  asymptotic distribution of the corner entry in
  rectangular tableaux} \label{sec:last-box}

Let $T$ be a Poissonized tableau of shape $\lambda$ and let the real number $u_0
    \in (\lambda_2-1, \lambda_1)$ be the $u$-coordinate of any point on the boundary
of the Young diagram in the first row. There is only one concave corner with
$u$-coordinate greater than $u_0$; this corner corresponds to the end of the
first row of $\lambda$. Schensted insertion $T \leftarrow z$ creates a new box
in this corner if and only if $z \geq T_{\lambda_1,1}$, i.e., if $z$ is larger
than  the last entry in the first row of $T$. It follows immediately that for
this choice of $u_0$, the value of the threshold is:
\[ F_T (u_0) = T_{\lambda_1,1}. \]
This coincides with the last entry in the first row of $T$.

Thanks to this observation, if $T$ is a uniformly random Poissonized tableau
with fixed shape $\lambda$, \cref{thm:wzormikolaja} provides convenient
information about the probability distribution of $T_{\lambda_1,1}$. In many
concrete cases, we may obtain interesting asymptotic results. We begin with the
following example.

\begin{corollary}[Asymptotic distribution of the corner entry in rectangular
        tableaux] \label{thm:last-box}

    Let $(p_l)$ and $(q_l)$ be sequences of positive integers such that $(p_l+ q_l)$
    tends to infinity and the limit \[ \alpha = \lim_{l \to \infty}
        \frac{q_l}{p_l+q_l } \] exists. We denote by \[ p_l\times q_l =
        (\underbrace{q_l,\dots,q_l}_{\text{$p_l$ times}})\] the rectangular Young
    diagram with $p_l$ rows and $q_l$ columns and by $n_l=p_l q_l$ the number of its
    boxes. Let $T^{(l)}$ be a uniformly random Poissonized tableau with shape
    $p_l\times q_l$. Define the normalized rightmost entry in the first row:
    \begin{equation}
        \label{eq:corner}
        Y^{(l)}:=\sqrt[4]{n_l} \left( T^{(l)}_{q_l,1}
        - \frac{q_l}{p_l+q_l} \right) \xrightarrow[l\to\infty]{d} N\left( 0,
        \sigma_\alpha  \right)
    \end{equation}
    converges to the centered Gaussian
    distribution with the variance \[ \sigma_\alpha^2=  \left[ \alpha (1-\alpha)
            \right]^{\frac{3}{2}}. \]
\end{corollary}
This result is due to Marchal
\cite{Marchal2016} who used very different methods. Below we present a new
proof.

\begin{proof}
    The diagram $p_l \times q_l$  exhibits two concave corners with
    $u$-coordinates $-p_l$ and $q_l$. Kerov's transition measure is supported at
    these corners, with respective probabilities $\frac{q_l}{p_l+q_l}$ and
    $\frac{p_l}{p_l+q_l}$.

    Let $u_0 \in (-p_l, q_l)$ lie between these concave corners. Applying
    \cref{thm:wzormikolaja}, we obtain the following expressions for the first two
    cumulants of $F_T(u_0)$:
    \begin{align}
        \E[F_T(u_0)]   & = \frac{q_l}{p_l+q_l},                                \\
        \Var[F_T(u_0)] & = \frac{q_l}{p_l+q_l} \cdot \frac{p_l}{p_l+q_l} \cdot
        \frac{1}{p_l+q_l+1}.
    \end{align}
    Consequently, after shift and scaling, the
    random variable $Y^{(l)}$ on the left-hand side of \eqref{eq:corner} has zero
    expectation, and its variance converges to $\sigma_\alpha^2$ as $l \to \infty$.

    Furthermore, \cref{thm:wzormikolaja} implies that:
    \[ \left| \kappa_k\left(
        F_T(u_0) \right) \right| \leq \frac{|\sT_{\nmoment}|}{(p_l+q_l+1)^{k-1}}. \]
    This leads to the following asymptotic behavior for the $k$-th cumulant of
    $Y^{(l)}$:
    \[ \kappa_k\left( Y^{(l)} \right) = O\left(
        (p_l+q_l)^{1-\frac{k}{2}} \right),\]
    which converges to zero for $k \geq 3$.
    Alternatively, the same conclusion can be obtained by applying
    \cref{coro:small-cumu}.

    \medskip

    We have thus demonstrated that the cumulants of $Y^{(l)}$ converge to their
    counterparts in the normal distribution on the right-hand side of
    \eqref{eq:corner}. This establishes convergence in moments. Since the normal
    distribution is uniquely determined by its moments, this moment convergence
    implies weak convergence of probability measures, proving the claim in
    \eqref{eq:corner}.
\end{proof}

The above method of proof is also applicable to the case when
\[
    \lambda^{(k)}=(\underbrace{q_{k,1},\dots,q_{k,1}}_{\text{$p_{k,1}$ times}},
    \underbrace{q_{k,2},\dots,q_{k,2}}_{\text{$p_{k,2}$ times}}, \dots,
    \underbrace{q_{k,i},\dots,q_{k,i}}_{\text{$p_{k,i}$ times}} ) \]
is a \emph{multi-rectangular Young diagram} obtained by stacking a fixed number
of rectangles.

It would be interesting to verify if \cref{thm:wzormikolaja}  could be applied
when $T$ is the random insertion tableau obtained by applying RSK to a sequence
of length $n$ of \iid random variables with the uniform distribution $U(0,1)$.
The goal would be to reprove the result of Azangulov \cite{AzangulovThesis2020}
that the random variable $n \left( 1- T_{\lambda_1,1} \right)$ converges to the
exponential distribution. See also \cite[Section~1.8]{MMS-Hammersley}.

\subfile{SECTION-trees-proof.tex}

\section{Acknowledgments}

Research was supported by Narodowe Centrum Nauki,
grant number \linebreak[4] 2017/26/A/ST1/00189. Additionally, the first named
author was supported by Narodowe Centrum Bada\'n i Rozwoju, grant number
POWR.03.05.00-00-Z302/17-00.

We thank %
Marek Bożejko, Maciej Dołęga, %
and Dan Romik for discussions and suggestions concerning the bibliography.

\printbibliography

\end{document}

%% file: figures/preamble-figures.tex
\colorlet{brazowy}{black}

\newcommand{\punkty}
{
	\coordinate (y00) at (0.4,2);
	\coordinate (y11) at (1.6,2);
	\coordinate (x1) at (1,0);
    \coordinate (x1) at (1,0);
    \coordinate (x2) at (2,0);
    \coordinate (x3) at (3,0);
    \coordinate (x4) at (4,0);
    \coordinate (x5) at (5,0);
    \coordinate (x6) at (6,0);
    \coordinate (x7) at (7,0);
    \coordinate (x8) at (8,0);
    \coordinate (x9) at (9,0);
    \coordinate (y1) at (1,2);
    \coordinate (y2) at (2,2);
    \coordinate (y23) at (3,1.4);
    \coordinate (y3) at (3,2);
    \coordinate (y4) at (4,2);
    \coordinate (y5) at (5,2);
    \coordinate (y6) at (6,2);
    \coordinate (yy2) at (2.3,2);
    \coordinate (yy3) at (3.7,2);
    \coordinate (y7) at (7,2);
    \coordinate (y8) at (8,2);
    \coordinate (y9) at (9,2);
    \coordinate (z1) at (1,4);
    \coordinate (z2) at (2,4);
    \coordinate (z23) at (3,2.6);
    \coordinate (z3) at (3,4);
    \coordinate (z4) at (4,4);
    \coordinate (gora) at (0,0.05);
    \coordinate (dol) at (0,-0.75);
    \coordinate (lewo) at (-.35,0.1);
    \coordinate (prawo) at (.35,0.1);
    \coordinate (lewodol) at (-.35,-0.7);
    \coordinate (prawodol) at (.35,-0.7);
}

\newcommand{\wierzcholek}[5]
{
    \draw[black,fill=#2] (#1) circle (4pt) node[above] at ($(#1)+(#4)$) {\tiny \large \textcolor{#5}{$#3$}};
}

\newcommand{\bwierzcholek}[5]
{
    \draw[black,fill=#2] (#1) circle (4pt) node[above] at ($(#1)+(#4)$) {\scriptsize   \textcolor{#5}{$#3$}};
}

\newcommand{\wierzcholekcbok}[5]
{
    \draw[red,fill=white] (#1) circle (4pt) node[above] at ($(#1)+(#4)$) {\tiny \scriptsize \textcolor{#5}{$#3$}};
   \prosto{$(#1)+(4pt ,0)$}{$(#1)+(- 4pt,0)$}{red}
   \prosto{$(#1)+(0, 4pt)$}{$(#1)+(0, - 4pt)$}{red}
}

\newcommand{\wierzcholekc}[5]
{
    \draw[red,fill=white] (#1) circle (4pt) node[above] at ($(#1)+(#4)$) {\tiny \large \textcolor{#5}{$#3$}};
   \prosto{$(#1)+(2.83 pt , 2.83 pt )$}{$(#1)+(- 2.83 pt, -2.83pt)$}{red}
   \prosto{$(#1)+(2.83 pt , -2.83 pt )$}{$(#1)+( -2.83 pt, 2.83pt )$}{red}
}

\newcommand{\wierzcholekpusty}[5]
{
    \draw[black,fill=#2] (#1) circle (0pt) node[above] at ($(#1)+(#4)$) {\tiny \large \textcolor{#5}{$#3$}};
}

\newcommand{\wierzcholekcpusty}[5]
{
    \draw[red,fill=#2] (#1) circle (0pt) node[above] at ($(#1)+(#4)$) {\tiny \large \textcolor{#5}{$#3$}};
    \draw[white,fill=white] (#1) circle (1.5pt) node[above] at ($(#1)+(#4)$) {};
}

\newcommand{\strzalka}[3]
{
    \draw[color=#3] [-Triangle] (#1) to [bend left=0] (#2);
}

\newcommand{\pstrzalka}[4]
{
    \draw[-Triangle, color=#4]  (#1) to [bend right=#3] (#2);
}

\newcommand{\lstrzalka}[4]
{
    \draw[-Triangle, color=#4] (#1) to [bend left=#3] (#2);
}

\newcommand{\prosto}[3]
{
    \draw[color=#3] (#1) to [bend left=0] (#2);
}

%% file: SECTION-teaser.tex
\section{Introduction}

\subsection{Basic definitions}

We begin by reviewing some fundamental combinatorial concepts. For a more
comprehensive treatment of this topic, we refer the reader to the book of Fulton
\cite{Fulton1997}.

\subsubsection{Young diagrams and tableaux}

A Young diagram is a finite collection of boxes arranged in the positive
quadrant, aligned to the left and bottom edges. This arrangement is known as the
French convention (see \cref{subfig:french}). To each Young diagram with $\ell$
rows, we associate an integer partition $\lambda = (\lambda_1, \ldots,
    \lambda_\ell)$, where $\lambda_j$ denotes the number of boxes in the $j$-th row,
counting from bottom to top. We identify a Young diagram with its corresponding
partition $\lambda$ and denote the total number of boxes by $|\lambda| =
    \lambda_1 + \cdots + \lambda_\ell$.

For asymptotic problems, it is convenient to draw Young diagrams using the
\emph{Russian convention} (see \cref{subfig:russian}). This corresponds to the
coordinate system $(u,v)$, which relates to the usual French Cartesian
coordinates as follows:
\begin{equation}
    \label{eq:Russian}
    u = x - y, \qquad v = x + y.
\end{equation}

\medskip

\begin{figure}[t]
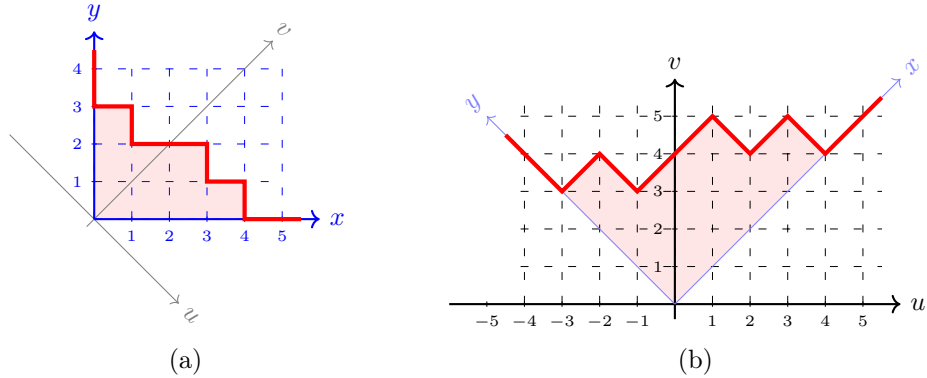

    \subfile{figures/FIGURE-diagrams.tex}
    \caption
    {   The Young diagram $(4,3,1)$ is depicted in two conventions:
        \protect\subref{subfig:french} the French convention,
        \protect\subref{subfig:russian} the Russian convention.
        In both representations, the solid red line illustrates the diagram's profile. 
        The coordinate systems are as follows:
        $(x,y)$ for the French convention, and
        $(u,v)$ for the Russian convention.
    }
    \label{fig:french}
\end{figure}

A \emph{tableau} (also known as \emph{semi-standard tableau}) is a filling of
the boxes of a Young diagram with numbers; we require that the entries should be
weakly increasing in each row (from left to right) and strictly increasing in
each column (from bottom to top). An example is given in \cref{fig:RSKa}. We say
that a tableau~$T$ of shape $\lambda$ is a \emph{standard Young tableau} if it
contains only entries from the set  $\{ 1,2,\dots,|\lambda|\}$ and each element
is used exactly once.

\subsubsection{The Schensted row insertion}

\begin{figure}[t]
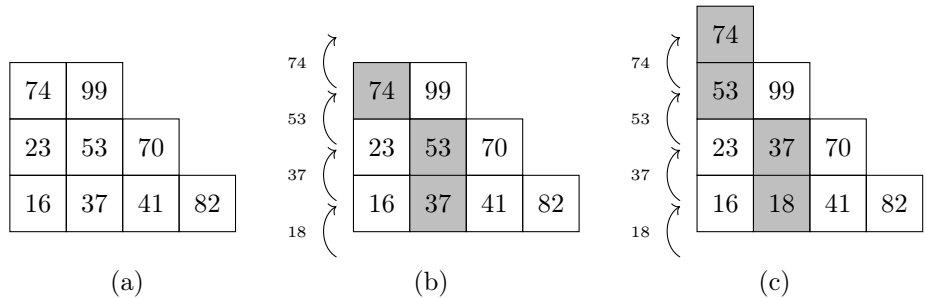

    \subfile{figures/FIGURE-RSK.tex}

    \caption{\protect\subref{fig:RSKa}~The original tableau $T$.
   	        \protect\subref{fig:RSKb}~The highlighted boxes indicate the bumping
				route for the Schensted insertion $T \leftarrow 18$. The numbers next to
				the arrows represent the bumped entries.
	        \protect\subref{fig:RSKc}~The resulting tableau after the Schensted
				insertion $T \leftarrow 18$.}

    \label{fig:RSK}
\end{figure}

The \emph{Schensted row insertion} is an algorithm that takes a tableau $T$ and
a number $z$ as input. The process begins by inserting $z$ into the first (bottom)
row of $T$, following these rules:
\begin{itemize}
    \item
          $z$ is placed in the leftmost box containing an entry strictly larger than $z$.

    \item
          If no such box exists, $z$ is appended to the end of the row in a new box, 
          and the algorithm terminates.

    \item
          If $z$ displaces an existing entry $z'$, this $z'$ is \emph{``bumped''} to the second row.

    \item
          The process repeats with $z'$ being inserted into the second row, following the same rules.

    \item
          This continues until a number is inserted into an empty box.
\end{itemize}
The resulting tableau is denoted as $T\leftarrow z$, see \cref{fig:RSKb,fig:RSKc} for an example.
The sequence of boxes whose contents change during this process is called \emph{the bumping route}.

Schensted insertion is a key component of the Robinson--Schensted--Knuth
algorithm (RSK), see below.

\subsubsection{The Robinson--Schensted--Knuth algorithm}
\label{sec:RSK-def}

This article considers a simplified version of the Robinson--Schensted--Knuth
algorithm (RSK), which is more accurately described as the Robinson--Schensted
algorithm. However, we retain the RSK acronym due to its widespread recognition.
The RSK algorithm maps a finite sequence $w=(w_1,\dots,w_{\nklatki})$ to a pair
of tableaux: the insertion tableau $P(w)$ and the recording tableau $Q(w)$.

The insertion tableau is defined as:
\begin{equation}
    \label{eq:insertion}
    P(w) = \Big( \big( (\emptyset \leftarrow w_1) \leftarrow
    w_2 \big) \leftarrow \cdots \Big) \leftarrow w_{\nklatki}.
\end{equation}
This represents the result of iteratively applying Schensted insertion to the
entries of $w$, beginning with an empty tableau $\emptyset$.

The recording tableau $Q(w)$ is a standard Young tableau with the same shape as
$P(w)$. Each entry in $Q(w)$ corresponds to the iteration number in
\eqref{eq:insertion} when that box was first filled. In other words, $Q(w)$
records the order in which the entries of the insertion tableau were populated.
Both $P(w)$ and $Q(w)$ share a common shape, denoted as $\RSK(w)$, which we
refer to as \emph{the RSK shape associated with $w$}.

\medskip

The RSK algorithm is a
fundamental tool in algebraic combinatorics and representation theory,
particularly in relation to Littlewood--Richardson coefficients (see
\cite{Fulton1997,Stanley1999}).

\subsection{The threshold for Schensted row insertion}

\subsubsection{The insertion function}

For a tableau $T$ and a real number $z$, we define
\[ \Ins(T; z) = (x,y) \]
as the
Cartesian French coordinates of the new box created by the Schensted row
insertion $T \leftarrow z$. This corresponds to the unique box in the skew
diagram:
\[ \operatorname{shape}(T \leftarrow z) / \operatorname{shape} T. \]

We define
\[ \uIns(T; z) = x-y\]
as the $u$-coordinate of $\Ins(T; z)$, as per equation \eqref{eq:Russian}.
For a fixed tableau $T$, we call the map
\[ z \mapsto \uIns(T; z) \]
\emph{the insertion function of $T$}. This function is weakly increasing.

\subsubsection{The threshold}

\begin{figure}
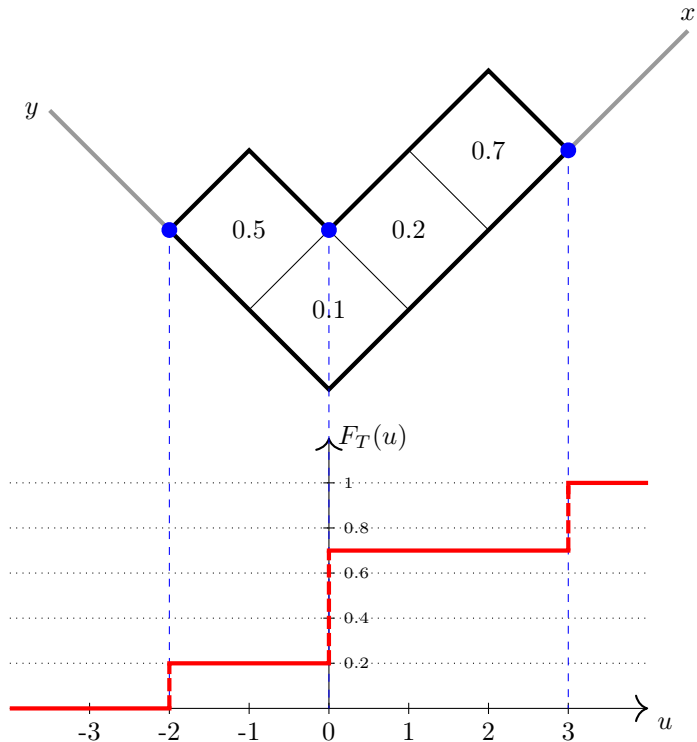

    \subfile{figures/FIGURE-cumulative-function-example.tex}

    \caption{A Poissonized tableau $T$ displayed in Russian coordinates. The
        threshold value $F_T(u)$, depicted by the red line along the bottom,
        illustrates its dependence on the $u$-coordinate.}
    \label{fig:threshold}
\end{figure}

We now focus on a Poissonized tableau $T$, where all entries belong to the unit
interval $[0,1]$. For $u_0 \in \R$, we define
\[ F_T(u_0) = \inf \big\{ z \in [0,1] : \uIns(T; z) > u_0 \big\} \]
as the threshold that separates small values $z$, for which the new box $\Ins(T; z)$ is
weakly to the left of $u_0$, from large values, for which the new box is
strictly to the right of $u_0$.
When the infimum is taken over an empty set, we define $F_T(u_0) = 1$ (see
\cref{fig:threshold} for an example).

\subsection{RSK algorithm and randomness}

\subsubsection{RSK applied to random input}

The investigation of the RSK algorithm applied to random input has proven to be
a fruitful area of research, revealing deep connections between combinatorics,
probability theory, and statistical mechanics. This approach has illuminated
links to the Plancherel growth process, directed last passage percolation
\cite{DauvergneNicaVirag}, and random polymer models
\cite{OconnellSeppalainenZygouras}.

A particularly significant achievement in this field is the solution to the
Ulam--Hammersley problem and its connection to the Tracy--Widom distribution from
random matrix theory  \cite{BaikDeift1999,Okounkov2000}. Romik's book
\cite{Romik2015} provides an excellent pedagogical introduction to these
concepts.

\subsubsection{Schensted row insertion into a random tableau}

Romik and \sniady \cite{RomikSniady2015,Sniady2014} pioneered a novel class of
problems involving the application of the RSK
algorithm to random data. Their work, rooted in ergodic theory, harmonic
analysis on the Young's graph, and the representation theory of the infinite
symmetric group, unexpectedly connected to the behavior of second-class
particles in interacting particle systems and the concept of competition
interfaces.

Using the notation established in this paper, these problems can be described as
follows: we begin with a random tableau $T$, drawn from a specified probability
distribution. We then insert a deterministic number $z$ into $T$ and examine the
position of the resulting new box, denoted as $\Ins(T;z)$. Our focus is on the
asymptotic behavior of this insertion process as the size of $T$ approaches
infinity.

\subsubsection{From first-order approximations to fluctuation analysis}

The aforementioned paper \cite{RomikSniady2015} provided first-order
approximations for the asymptotic behavior of $\Ins(T;z)$ as a function of $z$,
analogous to a law of large numbers. Our ultimate goal, to be achieved in a
forthcoming paper \cite{MarciniakSniady}, is to refine these results by
examining the fluctuations of $\Ins(T;z)$ around its mean value --- akin to a
central limit theorem. This current paper lays the technical foundation for
\cite{MarciniakSniady} by studying the probability distribution of the threshold
$F_T(u)$. By elucidating the finer probabilistic structure of random tableau
insertions, we open new avenues for understanding the interplay between
combinatorics, probability theory, and representation theory.

\subsection{Random Poissonized tableaux}
\label{sec:poissonized}

A \emph{Poissonized tableau} \cite{GorinRahman2019} is a tableau with entries
from the unit interval $[0,1]$. We denote the set of Poissonized tableaux with
shape $\lambda$ by $\mathcal{T}^{\lambda}$.

By numbering the boxes of $\lambda$ arbitrarily, each element of
$\mathcal{T}^{\lambda}$ can be identified with a point in the unit cube
$[0,1]^{\nklatki}$, where $n$ is the number of boxes in $\lambda$. The
requirement for increasing rows and columns corresponds to a set of inequalities
between coordinates, making $\mathcal{T}^{\lambda}$ identifiable with a convex
polytope in $[0,1]^{\nklatki}$. This polytope has positive volume, allowing us
to equip it with \emph{the uniform probability measure}. Thus, it makes sense to
speak about a \emph{uniformly random Poissonized tableau with shape $\lambda$}.

\subsection{Problem statement and paper overview}

This paper addresses a fundamental question in the study of random tableaux:

\begin{problem}[Distribution of the threshold]
\label{prob:main-problem}
Let $T$ be a uniformly random Poissonized tableau of
fixed shape $\lambda$, and let $u_0$ be a fixed real number. Characterize the
probability distribution of the random variable $F_T(u_0)$ in a manner that
facilitates asymptotic analysis as the size of the Young diagram $\lambda$
approaches infinity.
\end{problem}

Our main result, \cref{thm:wzormikolaja}, provides a comprehensive solution to
this problem. The paper is structured as follows:
\begin{itemize}
    \item \cref{sec:additional-definitions} introduces key concepts
          essential to our analysis, including Kerov's transition measure of a Young
          diagram and non-crossing alternating trees.

    \item \cref{sec:main-result} presents our main theorem (Theorem
          \ref{thm:wzormikolaja}), which provides an explicit combinatorial formula for
          the cumulants of the random variable $F_T(u_0)$.

    \item \cref{sec:last-box} demonstrates an application of our main
          result, examining the fluctuations of the last box in the first row of a
          uniformly random Poissonized tableau with a prescribed large shape.

    \item \crefrange{sec:dowod-kumulant-zaczyna-sie}{sec:dowod-kumulant-konczy-sie}
          contain the proof of the main result.

\end{itemize}

%% file: figures/FIGURE-diagrams.tex
	\subfloat[]{
		\begin{tikzpicture}
			
			\begin{scope}[scale=0.5/sqrt(2),rotate=-45,draw=gray]
				
				\begin{scope}[draw=gray,rotate=45,scale=sqrt(2)]
					\fill[fill=red!10] (4,0) -- (4,1) -- (3,1) -- (3,2) -- (1,2)
					-- (1,3) -- (0,3) -- (0,0) -- cycle ;
				\end{scope}
				
				\begin{scope}[rotate=45,draw=blue,scale=sqrt(2)]
					\draw[ultra thin, loosely dashed] (0,0) grid (5,4);
				\end{scope}
				
				\draw[->,thin] (-4.5,0) -- (4.5,0)
				node[anchor=west,rotate=-45]{\textcolor{gray}{$u$}};
				
				\draw[->,thin] (0,-0.4) -- (0,9.5)
				node[anchor=south,rotate=-45]{\textcolor{gray}{$v$}};
				
				\begin{scope}[draw=blue,rotate=45,scale=sqrt(2)]
					
					\draw[->,thick] (0,0) -- (6,0) node[anchor=west]{\textcolor{blue}{$x$}};
					\foreach \x in {1, 2, 3, 4, 5}
					{ \draw (\x, -2pt) node[anchor=north] {\textcolor{blue}{\tiny{$\x$}}} -- (\x,
						2pt); }
					
					\draw[->,thick] (0,0) -- (0,5) node[anchor=south] {\textcolor{blue}{$y$}};
					\foreach \y in {1, 2, 3, 4}
					{ \draw (-2pt,\y) node[anchor=east] {\textcolor{blue}{\tiny{$\y$}}} -- (2pt,\y); }
					
					\draw[ultra thick,draw=red] (5.5,0) -- (4,0) -- (4,1) -- (3,1) --
					(3,2) -- (1,2) -- (1,3) -- (0,3) -- (0,4.5) ;
					
				\end{scope}
				
			\end{scope}
		\end{tikzpicture}
		\label{subfig:french}
	}
	\hfill
	\subfloat[]
	{
		\begin{tikzpicture}
			\begin{scope}[xshift=7cm, yshift=-0.5cm, scale=0.5]
				
				\begin{scope}[draw=gray,rotate=45,scale=sqrt(2)]
					\fill[fill=red!10] (4,0) -- (4,1) -- (3,1) -- (3,2) -- (1,2)
					-- (1,3) -- (0,3) -- (0,0) -- cycle ;
				\end{scope}
				
				\begin{scope}
					\clip (-4.5,0) rectangle (5.5,5.5);
					\draw[ultra thin, loosely dashed] (-6,0.01) grid (6,6);
				\end{scope}
				
				\draw[->,thick] (-6,0) -- (6,0) node[anchor=west]{$u$};
				\foreach \z in {-5, -4, -3, -2, -1, 1, 2, 3, 4, 5}
				{ \draw (\z, -2pt) node[anchor=north] {\tiny{$\z$}} -- (\z, 2pt); }
				
				\draw[->,thick] (0,-0.4) -- (0,6) node[anchor=south]{$v$};
				\foreach \t in {1, 2, 3, 4, 5}
				{ \draw (-2pt,\t) node[anchor=east] {\tiny{$\t$}} -- (2pt,\t); }

				\begin{scope}[draw=blue!50,rotate=45,scale=sqrt(2)]
					
					\draw[->,thin] (0,0) -- (6,0) node[anchor=west,rotate=45]
					{\textcolor{blue!50}{{$x$}}};
					
					\draw[->,thin] (0,0) -- (0,5) node[anchor=south,rotate=45]
					{\textcolor{blue!50}{{$y$}}};
					
					\draw[ultra thick,draw=red] (5.5,0) -- (4,0) -- (4,1) -- (3,1) --
					(3,2) -- (1,2) -- (1,3) -- (0,3) -- (0,4.5) ;
					
				\end{scope}
			\end{scope}
			
		\end{tikzpicture}
		\label{subfig:russian}
	}

%% file: figures/FIGURE-RSK.tex
	\centering
	\hfill
	\subfloat[]
	{
		\begin{tikzpicture}[scale=0.75]
			\clip (-0.1,-0.5) rectangle (4.1,4.5);
			\draw (0,0) rectangle +(1,1); 
			\node at (0.5,0.5) {16}; 
			\draw (1,0) rectangle +(1,1); 
			\node at (1.5,0.5) {37}; 
			\draw (2,0) rectangle +(1,1); 
			\node at (2.5,0.5) {41}; 
			\draw (3,0) rectangle +(1,1); 
			\node at (3.5,0.5) {82}; 
			\draw (0,1) rectangle +(1,1); 
			\node at (0.5,1.5) {23}; 
			\draw (1,1) rectangle +(1,1); 
			\node at (1.5,1.5) {53}; 
			\draw (2,1) rectangle +(1,1); 
			\node at (2.5,1.5) {70}; 
			\draw (0,2) rectangle +(1,1); 
			\node at (0.5,2.5) {74}; 
			\draw (1,2) rectangle +(1,1); 
			\node at (1.5,2.5) {99}; 	
		\end{tikzpicture}
		\label{fig:RSKa}
	}
	\hfill
	\subfloat[]
	{
		\begin{tikzpicture}[scale=0.75]
			\clip (-1.5,-0.5) rectangle (4.1,4.5);
			\fill[blue!10] (1,0) rectangle +(1,1);
			\fill[blue!10] (1,1) rectangle +(1,1);
			\fill[blue!10] (0,2) rectangle +(1,1);
			\draw (0,0) rectangle +(1,1); 
			\node at (0.5,0.5) {16}; 
			\draw[fill=lightgray] (1,0) rectangle +(1,1); 
			\node at (1.5,0.5) {37}; 
			\draw (2,0) rectangle +(1,1); 
			\node at (2.5,0.5) {41}; 
			\draw (3,0) rectangle +(1,1); 
			\node at (3.5,0.5) {82}; 
			\draw (0,1) rectangle +(1,1); 
			\node at (0.5,1.5) {23}; 
			\draw[fill=lightgray] (1,1) rectangle +(1,1); 
			\node at (1.5,1.5) {53}; 
			\draw (2,1) rectangle +(1,1); 
			\node at (2.5,1.5) {70}; 
			\draw[fill=lightgray] (0,2) rectangle +(1,1); 
			\node at (0.5,2.5) {74}; 
			\draw (1,2) rectangle +(1,1); 
			\node at (1.5,2.5) {99}; 
			\draw[->] (-0.3,-0.45) to[bend left=60] (-0.3,0.45);
			\draw[->] (0,1) +(-0.3,-0.45) to[bend left=60] +(-0.3,0.45);
			\draw[->] (0,2) +(-0.3,-0.45) to[bend left=60] +(-0.3,0.45);
			\draw[->] (0,3) +(-0.3,-0.45) to[bend left=60] +(-0.3,0.45);
			\tiny
			\node[] at (-1,0) {18};
			\node[] at (-1,1) {37};
			\node[] at (-1,2) {53};
			\node[] at (-1,3) {74};
		\end{tikzpicture}
		\label{fig:RSKb}
	}
	\hfill
	\subfloat[]
	{
		\begin{tikzpicture}[scale=0.75]
			\clip (-1.5,-0.5) rectangle (4.1,4.5);
			\fill[blue!10] (1,0) rectangle +(1,1);
			\fill[blue!10] (1,1) rectangle +(1,1);
			\fill[blue!10] (0,2) rectangle +(1,1);
			\fill[blue!10] (0,3) rectangle +(1,1);
			\draw (0,0) rectangle +(1,1); 
			\node at (0.5,0.5) {16}; 
			\draw[fill=lightgray] (1,0) rectangle +(1,1); 
			\node at (1.5,0.5) {18}; 
			\draw (2,0) rectangle +(1,1); 
			\node at (2.5,0.5) {41}; 
			\draw (3,0) rectangle +(1,1); 
			\node at (3.5,0.5) {82}; 
			\draw (0,1) rectangle +(1,1); 
			\node at (0.5,1.5) {23}; 
			\draw[fill=lightgray] (1,1) rectangle +(1,1); 
			\node at (1.5,1.5) {37}; 
			\draw (2,1) rectangle +(1,1); 
			\node at (2.5,1.5) {70}; 
			\draw[fill=lightgray] (0,2) rectangle +(1,1); 
			\node at (0.5,2.5) {53}; 
			\draw (1,2) rectangle +(1,1); 
			\node at (1.5,2.5) {99}; 
			\draw[fill=lightgray] (0,3) rectangle +(1,1); 
			\node at (0.5,3.5) {74}; 
			\draw[->] (-0.3,-0.45) to[bend left=60] (-0.3,0.45);
			\draw[->] (0,1) +(-0.3,-0.45) to[bend left=60] +(-0.3,0.45);
			\draw[->] (0,2) +(-0.3,-0.45) to[bend left=60] +(-0.3,0.45);
			\draw[->] (0,3) +(-0.3,-0.45) to[bend left=60] +(-0.3,0.45);
			\tiny
			\node[] at (-1,0) {18};
			\node[] at (-1,1) {37};
			\node[] at (-1,2) {53};
			\node[] at (-1,3) {74};
		\end{tikzpicture}
		\label{fig:RSKc}
	}

%% file: figures/FIGURE-cumulative-function-example.tex
	\begin{tikzpicture}[scale=1.5,rotate=45]

		\coordinate (start) at (-2,-2);
		\coordinate (p1) at ($(start)+(-2,2)$);
		\coordinate (p2) at ($(start)+(2,-2)$);
		
		\coordinate (a) at ($2*(0.70710678118,0.70710678118)$);
		\coordinate (q0) at ($(p1)!(-4,0)!(p2)$);
		\coordinate (q1) at ($(p1)!(-2,0)!(p2)$);
		\coordinate (q2) at ($(p1)!(0,0)!(p2)$);
		\coordinate (q3) at ($(p1)!(3,0)!(p2)$);
		\coordinate (q4) at ($(p1)!(4,0)!(p2)$);

		\draw[decoration={markings,mark=at position 1 with {\arrow[scale=2]{>}}},
		postaction={decorate}] (p1) -- (p2) node[anchor=north west]{$u$};
		
		\draw[decoration={markings,mark=at position 1 with {\arrow[scale=2]{>}}},
		postaction={decorate}] (start) -- ($(start)+1.2*(a)$) node[anchor=west]{$F_T(u)$};

		\foreach \x in {-3, -2, -1, 0, 1, 2, 3} 
		{ \draw ($(p1)!(\x,0)!(p2)$) +(1pt,1pt) -- +(-1pt,-1pt) node[anchor=north] {\x}; }
		
		\foreach \y/\yt in {0.2/0.2,0.4/0.4,0.6/0.6,.8/0.8,1/1} 
		{ \draw[thin,dotted] ($(start)+\y*(a)+(-2,2)$)--($(start)+\y*(a)+(2,-2)$);
			\draw ($(start)+\y*(a)+(-1pt,1pt)$) -- ($(start)+\y*(a)+(1pt,-1pt)$) node [anchor=west]{\tiny$\yt$};
		}

		\draw[black!40, ultra thick] (4.5,0) node [anchor=south]{\textcolor{black}{$x$}} -- (0,0) -- (0,3.5) node[anchor=east]{\textcolor{black}{$y$}};

		\begin{scope}
			\clip (0,0) -- (3,0) -- (3,1) -- (1,1) -- (1,2) -- (0,2);
			\draw (0,0) grid (3,3); 
		\end{scope}

		\draw[ultra thick] (0,0) -- (3,0) -- (3,1) -- (1,1) -- (1,2) -- (0,2) -- cycle;

		\foreach \x/ \y/ \t in { 3/0/{(0, 0.2) }, 1/1/{(0.2,0.7)}, 0/2/{(0.7,1)} }
		{ 
			\fill[blue] (\x,\y) circle (2pt);
			\draw[blue, dashed] (\x,\y) -- ($(p1)!(\x,\y)!(p2)$ );
		}

		\draw[red,ultra thick] ($(q0)+0*(a)$)  -- ($(q1)+0*(a)$) -- ($(q1)+.01*(a)$); 
		\draw[red,ultra thick] ($(q1)+0.2*(a)-.01*(a)$) -- ($(q1)+0.2*(a)$) -- ($(q2)+0.2*(a)$) -- ($(q2)+0.2*(a)+.01*(a)$); 
		\draw[red,ultra thick] ($(q2)+0.7*(a)-.01*(a)$)  --($(q2)+0.7*(a)$)  -- ($(q3)+0.7*(a)$)-- ($(q3)+0.7*(a)+.01*(a)$); 
		\draw[red,ultra thick] ($(q3)+1*(a)-.05*(a)$)   --($(q3)+1*(a)$)   -- ($(q4)+1*(a)$);

		\draw[red,very densely dashed, ultra thick]  ($(q1)+0*(a)$) -- ($(q1)+0.2*(a)$); 
		\draw[red,very densely dashed, ultra thick] ($(q2)+0.2*(a)$) -- ($(q2)+0.7*(a)$); 
		\draw[red,very densely dashed, ultra thick] ($(q3)+0.7*(a)$) -- ($(q3)+1*(a)$);

		\path node () at (0.5,0.5){0.1};
		\path node () at (1.5,0.5){0.2};
		\path node () at (2.5,0.5){0.7};
		\path node () at (0.5,1.5){0.5};

	\end{tikzpicture}

%% file: SECTION-general-form.tex
\section{Additional definitions}
\label{sec:additional-definitions}

\subsection{Plancherel growth process}
\label{sec:plan}

Let $w_1, w_2, \dots$ be a sequence of \iid random variables with the uniform
distribution $U(0,1)$ on the unit interval $[0,1]$. Define
\[ \lambda^{(\nklatki)} := \RSK(w_1, \dots, w_{\nklatki}); \]
we refer to the random sequence of Young diagrams
\begin{equation}
    \label{eq:Plancherel-growth-process}
    \emptyset = \lambda^{(0)} \nearrow \lambda^{(1)} \nearrow \cdots
\end{equation}
as the \emph{Plancherel growth process} \cite[Chapter~1.19]{Romik2015}. It turns
out that \eqref{eq:Plancherel-growth-process} is a Markov chain; we will
describe its transition probabilities below.

\subsection{Transition measure of a Young diagram}
\label{sec:transition-measure}

\begin{figure}
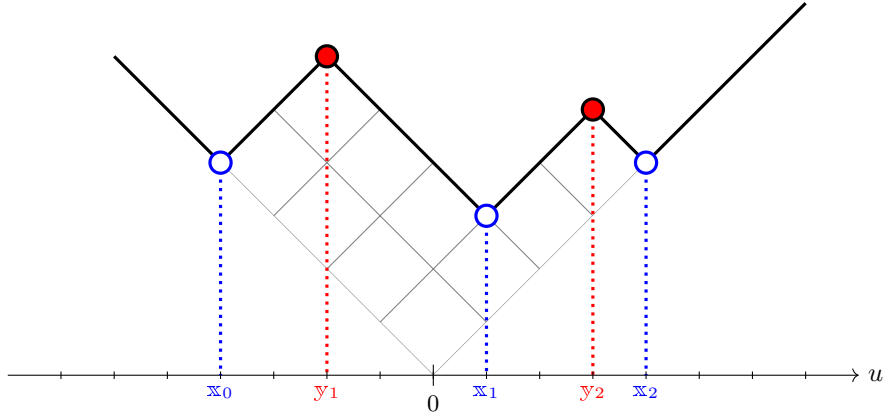

    \subfile{figures/FIGURE-cauchy.tex}
    \caption{Young diagram $(4,2,2,2)$ illustrating
        concave corners (unfilled circles) and
        convex corners (filled circles) with their $u$-coordinates.}
    \label{fig:cauchy}
\end{figure}

For a Young diagram $\lambda$ with $n$ boxes, let
$\kerx_0<\cdots<\kerx_{\kerell}$ denote the $u$-coordinates of its concave
corners and $\kery_1<\cdots<\kery_\kerell$ the $u$-coordinates of its convex
corners (see \cref{fig:cauchy}). The \emph{Cauchy transform} of $\lambda$ is
defined as the rational function \cite{Kerov1993,KerovBook}:
\begin{equation}
    \label{eq:cauchy-def}
    \cauchy_\lambda(z)=
    \frac{(z-\kery_1)\cdots(z-\kery_\kerell)}{(z-\kerx_0)\cdots(z-\kerx_\kerell)}.
\end{equation}
Note that in the work of Kerov this function is called \emph{the generating function} of $\lambda$.

The Cauchy transform can be uniquely expressed as a sum of simple fractions:
\[ \cauchy_\lambda(z)=\sum_{0\leq i\leq \kerell}\frac{p_i}{z-\kerx_i} \]
where $p_0,\dots,p_\kerell>0$ and $p_0+\cdots+p_\kerell=1$.
The \emph{transition measure} of $\lambda$ is defined as the discrete measure:
\[ \mu_\lambda=p_0\delta_{\kerx_0}+\cdots +p_\kerell\delta_{\kerx_\kerell} \]
such that $\cauchy_\lambda(z)$ is indeed its Cauchy transform:
\[ \cauchy_\lambda(z)=\int_{\R}\frac{1}{z-x}\dif \mu_\lambda(x). \]

Kerov showed that the transition probabilities of the Markov chain
\eqref{eq:Plancherel-growth-process} are encoded by this transition measure
\cite{Kerov1993,KerovBook}. Specifically, the conditional probability that the
new box will have $u$-coordinate $\kerx_i$ is given by:
\begin{equation}
    \label{eq:residuum}
    \Pro\giventhat*{ u\left(\lambda^{(n+1)} / \lambda^{(n)}\right)
        = \kerx_i }{ \lambda^{(n)}=\lambda }
    = p_i = \operatorname{Res}_{\kerx_i} \cauchy_\lambda.
\end{equation}
This probability corresponds to both the atom of the transition measure and the
residue of the Cauchy transform at $\kerx_i$.

%% file: figures/FIGURE-cauchy.tex
	
	\begin{tikzpicture}[scale=1,rotate=45]
		\coordinate (p0) at (4,0);
		\coordinate (p1) at (4,1);
		\coordinate (p2) at (2,1);
		\coordinate (p3) at (2,4);
		\coordinate (p4) at (0,4);
		\coordinate (p5) at (0,6);
		\begin{scope}
			\clip  (p0) -- (p1) -- (p2) -- (p3) -- (p4) -- (0,0);
			\draw[black!50] (0,0) grid (10,10);
		\end{scope}
		\draw[very thick] (7,0) -- (p0) -- (p1) -- (p2) -- (p3) -- (p4) -- (p5);
		\draw[->] (-4,4) -- (4,-4) node[anchor=west] {$u$};
		\draw[dotted,blue,very thick] (p0) -- (4/2,-4/2) node[anchor=north] {\small $\kerx_2$};
		\draw[dotted,blue,very thick] (p2) -- (1/2,-1/2) node[anchor=north] {\small $\kerx_1$};
		\draw[dotted,blue,very thick] (p4) -- (-4/2,4/2) node[anchor=north] {\small $\kerx_0$};
		\draw[dotted,red,very thick] (p1) -- (3/2,-3/2) node[anchor=north] {\small $\kery_2$};
		\draw[dotted,red,very thick] (p3) -- (-2/2,2/2) node[anchor=north] {\small $\kery_1$};
        \draw (0,0) +(0.1,0.1) -- +(-0.1,-0.1) node[anchor=north] {\small $0$};
		\draw[blue,very thick,fill=white] 
		(p0) circle (4pt)
		(p2) circle (4pt)
		(p4) circle (4pt);
		\draw[very thick,fill=red] 
		(p1) circle (4pt)
		(p3) circle (4pt);
		
		\foreach \u in {-7,...,7}
		{\draw (\u/2,-\u/2)  +(-1pt,-1pt) -- +(1pt,1pt); }
	\end{tikzpicture}

%% file: SECTION-cumulative-function.tex
\subsection{Cumulants and moments}
\label{sec:cumulants-and-moments}

Let $X$ be a random variable with the sequence of moments
\((m_{n})_{n=1}^{\infty}\), where \(m_{n} = \mathbb{E}[X^{n}]\). The formal
power series
\[
    \mathbb{E}[e^{tX}] = \sum_{n=0}^{\infty} \frac{m_{n}}{n!} t^{n}
\]
is its exponential moment generating function or \emph{formal Fourier--Laplace
	transform}. The coefficients \((\kappa_{n})_{n=1}^{\infty}\) of its formal
logarithm
\[
    \log \mathbb{E}[e^{tX}] = \sum_{n=1}^{\infty} \kappa_{n} \frac{t^{n}}{n!}
\]
are called the \emph{cumulants} \cite{Lauritzen2002} of the random variable X.
The first cumulant is the expected value and the second cumulant is the
variance:
\begin{align*}
    \kappa_1 & = \mathbb{E}[X],   \\
    \kappa_2 & = \mathrm{Var}(X).
\end{align*}

Cumulants are related to the moments via \emph{the moment-cumulant formula}
\begin{equation}
    \label{eq:rekurencja}
    m_{n} = \sum_{\pi \in \Pi_{n}} \prod_{b \in \pi} \kappa_{|b|},
\end{equation}
where the sum runs over all set-partitions \(\pi\) of the set \(\{1, \dots,
n\}\), and the product runs over all blocks of the partition \(\pi\). For
example, for \(n=3\), there are $5$ set-partitions of the set \(\{1, 2, 3\}\),
namely,
\[
    \left\{ \{1\}, \{2\}, \{3\} \right\}, \quad
    \left\{ \{1, 2\}, \{3\} \right\}, \quad
    \left\{ \{1, 3\}, \{2\} \right\}, \quad
    \left\{ \{2, 3\}, \{1\} \right\}, \quad
    \left\{ \{1, 2, 3\} \right\}.
\]
Therefore,
\[
    m_3 = \kappa_1^3 + \kappa_2 \kappa_1 + \kappa_2 \kappa_1 + \kappa_2 \kappa_1 + \kappa_3 = 
    \kappa_1^3 + 3 \kappa_2 \kappa_1 + \kappa_3.
\]

Cumulants are useful in probability theory because they provide a compact and
insightful way to characterize probability distributions. They transform simply
under affine transformations and can be used to quantify deviations from the
Gaussian law.

%% file: SECTION-trees.tex
\subsection{Directed, weighted graphs}

A \emph{directed graph} is a graph in which every edge has a direction. We
denote an edge from vertex $a$ to vertex $b$ as $(a,b)$. The directed graphs we
consider do not have multiple edges but may contain \emph{loops}, i.e., edges
of the form $(a,a)$. We assume that if $a \neq b$ and $(a,b)$ is an edge, then
the opposite edge $(b,a)$ is \emph{not} present.

The vertices in our graphs are
colored black, red, or white. For readers of the non-colored printed version,
red vertices are depicted as crossed-out circles (see \cref{fig:dwadwa}).

For a
graph $\g$, we use the following notation:
\begin{itemize}
    \item $\sV_{\g}$:
          set of all vertices,
    \item $\sB_{\g}$: set of black vertices,
    \item $\sR_{\g}$:
          set of red vertices,
    \item $\sW_{\g}$: set of white vertices,
    \item $\sE_{\g}$:
          set of edges.
\end{itemize}

A \emph{weighted graph} is a graph where each edge $e$ is assigned a numerical
\emph{weight}, denoted by $\weight(e) \in \R$.

\subsection{Decorations} Let $\Deco \subset \R$ be a fixed discrete set, and
let $u_0 \in \R$ be a fixed real number. Elements of the interval $(-\infty,
    u_0]$ are called \emph{small}, while elements of the interval $(u_0, \infty)$
are called \emph{big}.

For a given graph $\g$, a function $\mathbf{x} \colon V_{\g} \rightarrow \Deco$
is called \emph{a $u_0$-decoration of the graph $\g$} if the following
conditions are satisfied:
\begin{align*}
    \mathbf{x}(b) & \text{ is small for
    each } b \in \sB_{\g},              \\ \mathbf{x}(w) & \text{ is big for each } w \in
       \sW_{\g}.
\end{align*}
For simplicity, we denote $x_v = \mathbf{x}(v)$ for $v
    \in V_{\g}$. The set of all $u_0$-decorations of the graph $\g$ is denoted by
$D_{\g}(u_0)$. When the value of $u_0$ is clear from the context, we will
simply refer to \emph{decorations} and write $D_{\g} = D_{\g}(u_0)$.

\subsection{Non-crossing alternating trees}

\begin{figure}
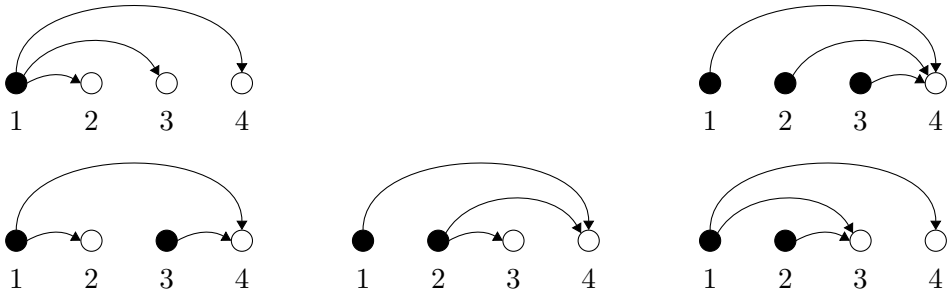

    \subfile{figures/FIGURE-tree-four.tex} \caption{All non-crossing
        alternating trees with $4$ vertices.} \label{fig:treeswithfour}
\end{figure}

Let $\nmoment \geq 1$ be a natural number. We define a \emph{non-crossing
    alternating tree} with $\nmoment$ vertices numbered $1,\dots,\nmoment$ as a
tree satisfying the following conditions:
\begin{enumerate}[label=(\alph*)]
    \item Each vertex is colored either black or white.

    \item
          \label{item:non-crossing-b}
          For any edge connecting vertices $b$ and $w$ with
          $b < w$, vertex $b$ is black and vertex $w$ is white.

    \item
          \label{item:noncrossing}
          There do not exist four vertices $v_1 < v_2 < v_3 <
              v_4$ such that $v_1$ is connected to $v_3$, and $v_2$ is connected to $v_4$.
\end{enumerate}
See \cite[Exercise~6.19(p) and its solution]{Stanley1999}, as
well as \cite[Section~6]{GelfandGraev1997}. Condition \ref{item:noncrossing}
has a natural graphical interpretation: when drawing the vertices on the real
line and the edges as arcs above it, the edges do not cross, see \cref{fig:treeswithfour}.

For the special case $\nmoment=1$, we define the unique non-crossing
alternating tree as consisting of a single black vertex (see
\cref{fig:llterm}). With this convention, condition~\ref{item:non-crossing-b}
ensures that the vertex coloring can be uniquely determined from the edge
information.

We denote by $\sT_{\nmoment}$ the set of all non-crossing alternating trees
with $\nmoment$ vertices. For example, $|\sT_3|=2$ (see
\cref{fig:treeswiththree}), and $|\sT_4|=5$ (see \cref{fig:treeswithfour}).

\begin{figure}
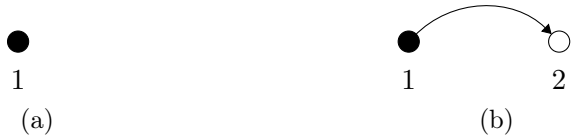

    \subfile{figures/FIGURE-tree-one-and-two.tex} \caption{
        \protect\subref{fig:llterm} The unique non-crossing alternating tree with one
        vertex. \protect\subref{fig:rrterm} The unique non-crossing alternating tree
        with two vertices. \label{fig:treeswithonetwo}
    }
\end{figure}

\begin{figure}
    \subfile{figures/FIGURE-tree-three.tex} \caption{All
        non-crossing alternating trees with $3$ vertices.} \label{fig:treeswiththree}
\end{figure}

\begin{remark}
    The number of non-crossing alternating trees with $n$ vertices
    is equal to the (shifted) Catalan number: \[ |\sT_n| = C_{n-1} =
        \frac{(2n-2)!}{n! (n-1)!} \] This demonstrates a connection to the rich
    combinatorics of Catalan objects \cite[Exercise 6.19, subpoint
        (p)]{Stanley1999}.
\end{remark}

In the following, we will treat any non-crossing alternating tree as a directed
graph, with each edge $(b,w)$ directed from the black vertex to the white
vertex, or equivalently, from the left vertex to the right vertex.

\section{The main result} \label{sec:main-result}

\subsection{Cumulants of the threshold}

Let $\lambda$ be a fixed Young diagram and $T$ be a uniformly random Poissonized
tableau of shape $\lambda$. Let $u_0\in\R$ be fixed.  We will now present a
closed formula for the cumulants of the threshold $F_T(u_0)$.

Assume that the discrete set $\Deco$ of decoration values contains the set of
$u$-coordinates of all concave corners of $\lambda$. In other words, we assume
that the support of the transition measure of $\lambda$ is contained in $\Deco$.
For example, we may take $\Deco = \Z$ to be the set of integers. Note that the
set $D_{\g}(u_0)$ used below depends implicitly on this choice of $\Deco$.

\begin{theorem}[The main result] \label{thm:wzormikolaja}
    With the above
    notations, for each $u_0 \in \R$ the $\nmoment$-th cumulant of the random
    variable $F_T(u_0)$ is given by
    \begin{equation}
        \label{eq:main-result-cumulants}
        \kappa_{\nmoment}\left(F_T\left(u_0\right)\right)=
         (\nmoment-1)!\sum_{\gT\in\sT_{\nmoment}} \sum_{\mathbf{x} \in D_{\gT}(u_0)} 
         \frac{(-1)^{|\sB_{\gT}|-1}\prod\limits_{j=1}^{\nmoment} \mu_\lambda(x_j) }{ \prod\limits_{(b,w)\in \sE_{\gT}}\left(x_w-x_b+w-b\right)},
    \end{equation}
    where $\mu_\lambda(x_j)$ denotes the probability corresponding to the atom $x_j$ of the transition measure~$\mu_\lambda$.
\end{theorem}

The complete proof spans Sections \ref{sec:dowod-kumulant-zaczyna-sie} through
\ref{sec:dowod-kumulant-konczy-sie}. In \cref{sec:ideas}, we present a concise
overview of the key ideas underlying the demonstration.

\begin{example}
    The first two cumulants, the expected value and the variance of $F_T(u_0)$ are given by
    \begin{align}
        \label{eq:expval}
        \E  F_T\left(u_0\right)  =    & \sum_{x_1 \leq
        u_0} \mu_\lambda(x_1),                         \\
        \label{eq:variance}
        \Var  F_T \left(u_0 \right) = &
        \sum_{\substack{ x_1 \leq u_0                  \\ x_2 > u_0 }}
        \frac{1}{x_2-x_1+1} \ \mu_\lambda(x_1)\ \mu_\lambda(x_2),
    \end{align}
    respectively. The unique
    summand on the right-hand side of \cref{eq:expval} corresponds to the tree %
    in \cref{fig:llterm}, and the unique summand on the right-hand side of
    \cref{eq:variance} corresponds to the tree %
    in \cref{fig:rrterm}. The third cumulant of $F_T(u_0)$ is given by
    \begin{multline*}
        k_3\left( F_T\left(u_0\right) \right) = \sum_{\substack{ x_1
                \leq u_0 \\ x_2,x_3 > u_0 }}  \frac{2}{(x_2-x_1+1)(x_3-x_1+2)}
        \mu_\lambda(x_1)\ \mu_\lambda(x_2) \ \mu_\lambda(x_3)  \\ - \sum_{\substack{
                x_1, x_2 \leq u_0 \\ x_3 > u_0 }}  \frac{2}{(x_3-x_2+1)(x_3-x_1+2)}
        \mu_\lambda(x_1)\ \mu_\lambda(x_2)\  \mu_\lambda(x_3),
    \end{multline*}
    where
    the first summand on the right-hand side corresponds to the tree in
    \cref{fig:lterm}, and the second summand corresponds to the tree in
    \cref{fig:rterm}.
\end{example}

\begin{remark}
    \label{rem:cumulant-as-expected-value}
    The right-hand side of
    \eqref{eq:main-result-cumulants} can be interpreted as the expected value of
    the random variable $Z$ defined in the following way. Let
    $x_1,\dots,x_{\nmoment}$ be a sequence of \iid random variables, with the
    distribution given by the transition measure~$\mu_\lambda$. Let
    $\sT_{\nmoment}^{\mathbf{x}}$ denote the set of all trees $\gT \in
        \sT_{\nmoment}$ such that $\mathbf{x}=(x_1,\dots,x_{\nmoment})$ is a
    $u_0$-decoration of the tree $\gT$, i.e., $\mathbf{x} \in D_{\gT}$. The
    aforementioned random variable is defined as
    \[  %
        Z=(\nmoment-1)! \sum_{ \gT \in \sT_{\nmoment}^{\mathbf{x}} }
        \frac{(-1)^{|\sB_{\gT}|-1}}{ \prod\limits_{(b,w)\in
                \sE_{\gT}}\left(x_w-x_b+w-b\right)}. \]
\end{remark}

\subsection{Rational functions associated to a graph}
\label{sec:rational}

For an directed weighted graph $\g$ with the vertex set $\sV_{\g}=\{ v_1, \dots,
    v_t \}$, we consider the rational function \[ f_{\g} =
    f_{\g}(x_{v_1},\dots,x_{v_{t}}) = \frac{1}{\prod\limits_{e=(i,j)\in \sE_{\g}}
        [x_j-x_i+\weight(e)]} \in\mathbb{Q}(x_{v_1},\dots,x_{v_{t}})	\] in the variables
corresponding to the vertices of $\g$. In the following we will usually consider
the special case when $\g$ has the vertex set $V_{\g} =\{1,\dots,\nmoment\}$ so
that \[ f_{\g} = f_{\g}(x_1,\dots,x_\nmoment) %
    \in\mathbb{Q}(x_1,\dots,x_{\nmoment}).\]

For each tree $\gT \in \sT_{k}$, we define the weight of an edge $e=(i,j) \in
    \sE_{\gT}$
\begin{equation}
    \label{eq:konwencja-wag}
    \weight(e)=\weight(i,j)=j-i
\end{equation}
as the difference of the endpoints. With this convention
\eqref{eq:main-result-cumulants} can be written more compactly as
\begin{equation}
    \label{eq:reformulated}
    \kappa_{\nmoment}\left(F_T\left(u_0\right)\right)=
    (\nmoment-1)!\sum_{\gT\in\sT_{\nmoment}} \sum_{\mathbf{x} \in D_{\gT}(u_0)}
    (-1)^{|\sB_{\gT}|-1} f_{\gT}(x_1,\dots,x_{\nmoment})
    \prod\limits_{j=1}^{\nmoment} \mu_\lambda(x_j).
\end{equation}

\subsection{Upper bound for the cumulants}

For some asymptotic problems, we do not need the full power of
\cref{thm:wzormikolaja}. Instead, we only require explicit formulas for the mean
value ($n=1$), given by \eqref{eq:expval}, and for the variance ($n=2$), given
by \eqref{eq:variance}. Additionally, a rough estimate on the higher-order
cumulants for $n \geq 3$ is sufficient. The following result provides such an
upper bound. Although we will not use this estimate in the current paper, it
will be crucial for the forthcoming paper \cite{MarciniakSniady}.

\begin{corollary}
    \label{coro:small-cumu}
    We maintain the notations from
    \cref{thm:wzormikolaja}. For any $k\geq 1$
    \[ \left| \kappa_{\nmoment} \left(
        F_T(u_0) \right) \right|  \leq (\nmoment-1)! \left[ \cauchy^+_{\lambda}(u_0)
        \right]^{k-1}, \]
    where $\cauchy^+_{\lambda}(u_0)$ denotes a variant of the
    Cauchy transform of $\mu_\lambda$, defined as:
    \[
        \cauchy^+_{\lambda}(u_0)=\sum_{z}\frac{1}{|u_0-z|+1}\mu_\lambda(z).\] This
    variant of the Cauchy transform employs a kernel that is a regularized version
    of the absolute value of the standard kernel.
\end{corollary}
\begin{proof}
    We
    bound each summand in \cref{thm:wzormikolaja} via \cref{lem:cumulants-small}.
\end{proof}

\begin{lemma} \label{lem:cumulants-small}

    Let \( G = (V, E) \) be a directed, weighted bipartite tree with \( k \geq 1 \)
    vertices, where:
    \begin{itemize}
        \item Each vertex is colored either black or
              white. \item Each edge connects vertices of opposite colors. \item Edges are
              directed from black vertices to white vertices. \item Edge weights are real
              numbers \(\geq 1\).
    \end{itemize}

    Let \(\mu\) be a discrete probability measure on \(\mathbb{R}\), and let \(u_0
    \in \mathbb{R}\).

    Define \( S(G) \) as: \[ S(G) := \sum_{\mathbf{x}} \frac{\prod_{v \in V}
        \mu(\mathbf{x}(v))}{\prod_{e=(i,j) \in E} (\mathbf{x}(j) - \mathbf{x}(i) +
        w(e))} \] where the sum runs over all \(u_0\)-decorations \(\mathbf{x}\) of the
    vertices of \(G\).

    Then, the following inequality holds: \[ S(G) \leq
        \left[\cauchy^+_\mu(u_0)\right]^{k-1}. \]
\end{lemma}
\begin{proof}
    We use
    induction over the number of the vertices.

    For $k=1$ there is nothing to prove.

    For $k\geq 2$ let $w$ be a leaf of $G$ and let $G'$ be the tree $G$ after
    removal of the vertex $w$ and the adjacent edge. By a straightforward bound on
    the factor which corresponds to the unique edge adjacent to $w$ it follows that
    \[ S(G) \leq S(G')\  \cauchy^+_{\mu}(u_0) \] and the inductive step follows
    immediately.
\end{proof}

\subsection{Sketch of the proof of \cref{thm:wzormikolaja}} \label{sec:ideas}

\subsubsection{Dual conditioning}
The Plancherel growth process \eqref{eq:Plancherel-growth-process} is a random
walk on the set of Young diagrams. It can be viewed as the result of applying
the RSK algorithm to a sequence of \iid random
variables $w_1, w_2, \dots$.

Let $n \geq 0$ and $k \geq 1$ be fixed integers.
We obtain the \emph{anti-Pieri growth process} by conditioning in one of two equivalent ways:
\begin{itemize}
    \item
          We require specific entries from the underlying sequence $(w_n)$ to form a decreasing
          sequence:
          \begin{equation}
              \label{eq:Pieri1}
              w_{n+1} > w_{n+2} > \cdots > w_{n+k}.
          \end{equation}

    \item Alternatively, we condition the Plancherel growth process and require that
          the new boxes created during $k$ transitions
          \begin{equation}
              \label{eq:Pieri2}
              \lambda^{(n)} \nearrow \cdots \nearrow \lambda^{(n+k)}
          \end{equation}
          satisfy a simple
          geometric condition: they should form a sequence of boxes on the plane with
          decreasing $u$-coordinates.
\end{itemize}
This type of growth process was introduced by Romik and \sniady
\cite[Section~4.1]{RomikSniady2015}. In fact, they considered a \emph{Pieri
    growth process} for which the inequalities in \eqref{eq:Pieri1} are in the
opposite order, and similarly, the $u$-coordinates of the new boxes in
\eqref{eq:Pieri2} should be increasing.

The dual nature of the anti-Pieri growth process allows us to translate
quantities related to the underlying sequence $(w_n)$ into the language of
transition probabilities of the (conditioned) Plancherel growth process.
\Cref{lem:string}\ref{item:A} exemplifies this translation, relating the moments
of the threshold $F_T(u_0)$ to probabilities associated with the Plancherel
growth process.

This concept of dual conditioning --- applied to random sequences on one side
and random growth processes (or random tableaux) on the other --- with the RSK
correspondence serving as a bridge between these two realms, played a crucial
role also in the research of Maślanka and \sniady \cite{MaslankaSniady}.

\subsubsection{Transition probabilities for next steps}

While the transition probabilities for the first step of the Plancherel growth
process starting from $\lambda$ are given by the transition measure
$\mu_\lambda$, subsequent steps are determined by the transition measure of the
evolved Young diagram. Despite this complexity, \Cref{lem:string}\ref{item:C}
demonstrates how these later transition probabilities relate to the Cauchy
transform $\cauchy_\lambda$ of the original diagram $\lambda$, albeit through
intricate rational functions.

\subsubsection{Translation from rational functions to graphs}

In \cref{sec:simple-fractions}, we develop a combinatorial framework that
associates rational functions with directed, weighted graphs. This approach
allows certain algebraic operations on rational functions to be interpreted
combinatorially as simple graph manipulations, such as edge removal. In this way
we manage to write the aforementioned rational function as a sum of simple
fractions which are indexed by particularly simple directed graphs.

\subsubsection{Regularization for uniform treatment}

The formulas we derive for the moments involve summation over numerous
combinatorial terms. For certain Young diagrams, some of these terms must be
explicitly excluded due to singularities. Treating these problematic terms
separately would significantly complicate the proof. In \cref{rozdzial6}, we
introduce regularization techniques that allow for uniform treatment of all
terms, thereby streamlining the analysis.

\subsubsection{Combinatorial tour de force}

Finally, in \cref{sec:dowod-kumulant-konczy-sie}, we complete the proof of
\cref{thm:wzormikolaja} by employing combinatorial techniques involving graphs.
Specifically, we leverage the fact that cumulants correspond to considering only
connected graphs.

\subsection{The key idea behind \cref{thm:wzormikolaja}}

It is remarkable that the cumulants of the threshold yield such an elegant
closed formula as provided by \cref{thm:wzormikolaja}. This elegance stems from
the deep combinatorial structures underlying the problem.

The proof of this formula relies on Kerov's product formula
\eqref{eq:cauchy-def} for the Cauchy transform of the transition measure of a
Young diagram. This product formula is a powerful tool for the harmonic analysis
on the Young graph and the random walks on the set of Young diagrams.

Ultimately, Kerov's product formula is derived from the Frame--Robinson--Thrall
formula (\emph{the hook length formula}) \cite{hook-length-formula}, which
calculates the number of standard Young tableaux of a given shape. The hook
length formula itself is a cornerstone in combinatorial mathematics, with
applications ranging from representation theory to probability and algorithm
analysis.

Therefore, the formula we prove is fundamentally a consequence of the hook
length formula, highlighting the deep interplay between combinatorial structures
and probabilistic measures.

%% file: figures/FIGURE-tree-one-and-two.tex
	    \begin{center}
		\subfloat[]
		{ \begin{tikzpicture}[scale=1.0]
				\punkty	
				\wierzcholekpusty{x1}{white}{}{dol}{white}
				\wierzcholek{x2}{black}{1}{dol}{black}
				\wierzcholekpusty{x3}{white}{}{dol}{white}
			\end{tikzpicture}
			\label{fig:llterm}
		}
		\qquad
		\qquad
		\qquad
		\qquad
		\qquad
		\subfloat[]
		{
			\begin{tikzpicture}[scale=1.0]
				\punkty	
				\lstrzalka{$(x1)+(2.12pt,2.12pt)$}{$(x3)+(-2.83pt,2.83pt)$}{45}{black}
				\wierzcholek{x1}{black}{1}{dol}{black}
				\wierzcholekpusty{x2}{white}{}{dol}{white}
				\wierzcholek{x3}{white}{2}{dol}{black}
			\end{tikzpicture}
			\label{fig:rrterm}
		}
	\end{center}
	

%% file: figures/FIGURE-tree-three.tex
	    \begin{center}
		\subfloat[]
		{
			\begin{tikzpicture}[scale=1.0]
				\punkty	
				\lstrzalka{$(x1)+(4pt,0)$}{$(x2)+(-4pt, 0)$}{30}{black}
				\lstrzalka{$(x1)+(2.83pt,2.83pt)$}{$(x3)+(-2.83pt,2.83pt)$}{60}{black}	
				\wierzcholek{x1}{black}{1}{dol}{black}
				\wierzcholek{x2}{white}{2}{dol}{black}
				\wierzcholek{x3}{white}{3}{dol}{black}
			\end{tikzpicture}
			\label{fig:lterm}
		}
		\qquad
		\subfloat[]
		{
			\begin{tikzpicture}[scale=1.0]
				\punkty	
				\lstrzalka{$(x2)+(4pt,0)$}{$(x3)+(-4pt, 0)$}{30}{black}
				\lstrzalka{$(x1)+(2.83pt,2.83pt)$}{$(x3)+(-2.83pt,2.83pt)$}{60}{black}
				\wierzcholek{x1}{black}{1}{dol}{black}
				\wierzcholek{x2}{black}{2}{dol}{black}
				\wierzcholek{x3}{white}{3}{dol}{black}
			\end{tikzpicture}
			\label{fig:rterm}
		}
	\end{center}
	

%% file: figures/FIGURE-tree-four.tex
	        \begin{tikzpicture}[scale=1.0]
		\punkty	
		\lstrzalka{$(x1)+(4pt,0)$}{$(x2)+(-4pt,0)$}{30}{black}
		\lstrzalka{$(x1)+(2.83 pt, 2.83 pt)$}{$(x3)+(-2.83 pt, 2.83 pt)$}{60}{black}
		\lstrzalka{$(x1)+(0,4pt)$}{$(x4)+(0,4pt)$}{90}{black}	
		\wierzcholek{x1}{black}{1}{dol}{black}
		\wierzcholek{x2}{white}{2}{dol}{black}
		\wierzcholek{x3}{white}{3}{dol}{black}
		\wierzcholek{x4}{white}{4}{dol}{black}
	\end{tikzpicture}
	\hfill
	\begin{tikzpicture}[scale=1.0]
		\punkty	
		\lstrzalka{$(x3)+(4pt,0)$}{$(x4)+(-4pt,0)$}{30}{black}
		\lstrzalka{$(x2)+(2.83 pt, 2.83 pt)$}{$(x4)+(-2.83 pt, 2.83 pt)$}{60}{black}
		\lstrzalka{$(x1)+(0,4pt)$}{$(x4)+(0,4pt)$}{90}{black}	
		\wierzcholek{x1}{black}{1}{dol}{black}
		\wierzcholek{x2}{black}{2}{dol}{black}
		\wierzcholek{x3}{black}{3}{dol}{black}
		\wierzcholek{x4}{white}{4}{dol}{black}
	\end{tikzpicture}
	\\
	\begin{tikzpicture}[scale=1.0]
		\punkty	
		\lstrzalka{$(x1)+(4pt,0)$}{$(x2)+(-4pt,0)$}{30}{black}
		\lstrzalka{$(x3)+(4pt,0)$}{$(x4)+(-4pt,0)$}{30}{black}
		\lstrzalka{$(x1)+(0,4pt)$}{$(x4)+(0,4pt)$}{90}{black}			
		\wierzcholek{x1}{black}{1}{dol}{black}
		\wierzcholek{x2}{white}{2}{dol}{black}
		\wierzcholek{x3}{black}{3}{dol}{black}
		\wierzcholek{x4}{white}{4}{dol}{black}
	\end{tikzpicture}
	\hfill
	\begin{tikzpicture}[scale=1.0]
		\punkty	
		\lstrzalka{$(x2)+(4pt,0)$}{$(x3)+(-4pt,0)$}{30}{black}
		\lstrzalka{$(x1)+(0,4pt)$}{$(x4)+(0,4pt)$}{90}{black}	
		\lstrzalka{$(x2)+(2.83 pt, 2.83 pt)$}{$(x4)+(-2.83 pt, 2.83 pt)$}{60}{black}	
		\wierzcholek{x1}{black}{1}{dol}{black}
		\wierzcholek{x2}{black}{2}{dol}{black}
		\wierzcholek{x3}{white}{3}{dol}{black}
		\wierzcholek{x4}{white}{4}{dol}{black}
	\end{tikzpicture}
	\hfill
	\begin{tikzpicture}[scale=1.0]
		\punkty	
		\lstrzalka{$(x2)+(4pt,0)$}{$(x3)+(-4pt,0)$}{30}{black}
		\lstrzalka{$(x1)+(0,4pt)$}{$(x4)+(0,4pt)$}{90}{black}	
		\lstrzalka{$(x1)+(2.83 pt, 2.83 pt)$}{$(x3)+(-2.83 pt, 2.83 pt)$}{60}{black}	
		\wierzcholek{x1}{black}{1}{dol}{black}
		\wierzcholek{x2}{black}{2}{dol}{black}
		\wierzcholek{x3}{white}{3}{dol}{black}
		\wierzcholek{x4}{white}{4}{dol}{black}
	\end{tikzpicture}
	

%% file: SECTION-trees-proof.tex
\section{Anti-Pieri growth process} \label{sec:dowod-kumulant-zaczyna-sie}

\subsection{RSK as a source of uniformly ransom Poissonized tableaux}
\label{sec:proof-conditioning}

\begin{lemma} \label{lem:conditioning}

    Let $w=(w_1,\dots,w_{\nklatki})$ be a sequence of \iid~random variables with
    the uniform $U(0,1)$ distribution and let $\lambda$ be a Young diagram with
    $\nklatki$ boxes.

    The \emph{conditional} probability distribution of the insertion tableau $P(w)$
    \emph{under the condition that $\RSK(w)=\lambda$} coincides with the uniform
    probability distribution on~$\mathcal{T}^\lambda$. \end{lemma}
\begin{proof}
    Let us consider the unit cube \[
        [0,1]^{\nklatki}=\{(w_1,\dots,w_{\nklatki}):w_1,\dots,w_{\nklatki}\in
        [0,1]\},\] equipped with the Lebesgue measure. We remove all hyperplanes
    defined by
    \[ w_i-w_j=0 \qquad \text{for } 1\leq i<j\leq \nklatki.\]
    This removal is
    inconsequential from a measure-theoretic perspective, as these hyperplanes have
    Lebesgue measure zero. These hyperplanes partition the cube into $\nklatki!$
    isometric simplices, each with volume $\frac{1}{\nklatki!}$. There exists a
    bijective correspondence between these simplices and permutations in
    $\Sy{\nklatki}$. Each simplex $\simplex_{\sigma}$ comprises vectors with a
    prescribed linear order among their coordinates.

    For a given Young diagram $\lambda$, let $\mathcal{T}^\lambda$ denote the set of
    Poissonized tableaux of shape $\lambda$, endowed with the Lebesgue measure. For
    simplicity, we exclude tableaux with repeated entries from this set, which is
    again measure-theoretically insignificant.

    Under these notations, the Robinson--Schensted correspondence establishes a
    bijection between the aforementioned cube $[0,1]^{\nklatki}$ and the disjoint
    sum
    \begin{equation}
        \label{eq:RSK-is-a-bijection}
        \bigsqcup_{\lambda \vdash n}
        \mathcal{T}^\lambda \times f^\lambda,
    \end{equation}
    where $f^\lambda$
    represents the set of standard Young tableaux of shape $\lambda$. The second
    component of this correspondence, the map $Q$, when restricted to the simplex
    $\simplex_\sigma$, is constant and equal to the recording tableau $Q(\sigma)$.
    On the other hand, the first component, the map $P$, when restricted to
    $\simplex_\sigma$, arranges the entries of $(w_1,\dots,w_{\nklatki})$ into the
    boxes of the diagram $\RSK(\sigma)$.

    Consequently, the Robinson--Schensted correspondence is a piecewise isometry. It
    follows that this correspondence is a measure-preserving map when we equip
    $f^\lambda$ with the counting measure and each summand in
    \eqref{eq:RSK-is-a-bijection} with the product measure.

    Conditioning on the event $$\RSK(w) = \lambda$$ is therefore equivalent to
    considering the uniform measure (or, alternatively, the product measure
    multiplied by the scalar factor $\frac{n!}{(f^\lambda)^2}$) on a specific
    summand of \eqref{eq:RSK-is-a-bijection}, namely $$ \mathcal{T}^\lambda \times
        f^\lambda, $$ which concludes our proof.
\end{proof}

\subsection{Anti-Pieri growth process}

We will use the following notation which is intended as an analogue of the
falling factorial \[ \cauchy_\lambda^{\underline{k}}(x)= \begin{cases}
        \underbrace{\cauchy_\lambda(x)\ \cauchy_\lambda(x-1) \ \cdots \
        \cauchy_\lambda(x-k+1)}_{\text{$k$ factors}} & \text{if } k\geq 1, \\ 1 &
           \text{if } k=0\end{cases} \] for an integer $k\geq 0$.

\begin{lemma} \label{lem:string}
    Let $\lambda$ be a fixed Young diagram with $\nklatki$ boxes and $k \geq 1$ be
    an integer. Consider the Plancherel growth process starting at $\lambda$:
    \[ \lambda =\diagram_{\nklatki}\nearrow \diagram_{\nklatki+1}\nearrow \cdots
        \nearrow \diagram_{\nklatki+k}.\]
    Let $\mathbf{U} = (U_1, \dots, U_k)$ be the
    sequence of $u$-coordinates of the boxes added in each step, i.e.,
    \[ U_i=u(\diagram_{\nklatki+i}\setminus \diagram_{\nklatki+i-1})\]
    for $i \in \{1, \dots, k\}$.
    \begin{enumerate}[label=(\alph*)]
        \item \label{item:A}
              Let $T$ be a random Poissonized tableau of shape
              $\lambda$. Then for any $u_0 \in \mathbb{R}$, the $k$-th moment of the random
              variable $F_T(u_0)$ is given by:
              \begin{align*}
                  m_k \left( F_T \left(u_0 \right) \right) & = \E \left[ \big(
                  F_T(u_0) \big)^k \right]                                                                \\
                                                           & = k!\ \Pro\big(u_0\geq U_1>\cdots>U_k\big) .
              \end{align*}

        \item \label{item:B}
              If $U_1 > \cdots > U_k$, then the tuple $\mathbf{U} = (U_1, \dots, U_k)$ can be uniquely written as:
              \[			\mathbf{x}^{\underline{\mathbf{a}}}=
              \big(\underbrace{x_1,x_1-1,\dots,x_1-\h_1+1}_{\text{$\h_1$ times}},\dots,
                   \underbrace{x_\ell,x_\ell-1,\dots,x_\ell-\h_\ell+1}_{\text{$\h_\ell$ times}}\big), \]
              where $\mathbf{x} = (x_1, \dots, x_\ell)$ and $\mathbf{a} = (\h_1, \dots,
                  \h_{\ell})$, with $x_1 > \cdots > x_\ell$ being $u$-coordinates of some
              concave corners of $\lambda$ and $\h_1, \dots, \h_\ell \geq 1$ integers
              such that $\h_1 + \cdots + \h_{\ell} = k$ (see \cref{fig:pieri}).

        \item
              \label{item:C}
              Let $x_1, \dots, x_\ell$ be the $u$-coordinates of some
              concave corners of $\lambda$, and let $\h_1, \dots, \h_\ell \geq 1$ be
              integers such that $\h_1 + \cdots + \h_{\ell} = k$. Assume the following
              condition holds:
              \begin{enumerate}[label=(X)]
                  \item \label{item:annoying-condition}
                        For each $i \in \{1, \dots, \ell\}$,
                        the set $\{x_i-1, x_i-2, \dots, x_i-a_i+1\}$ is disjoint from the set of
                        $u$-coordinates of the concave corners of $\lambda$.
              \end{enumerate}
              Then:
              \begin{equation}
                  \label{eq:big-product}
                  \mathbb{P}\Big[\mathbf{U} = \mathbf{x}^{\underline{\mathbf{a}}}\Big] =
                  \func(x_1, \dots, x_\ell)
                  \prod_{1 \leq i \leq \ell} \frac{(-1)^{\h_i-1}}{\h_i} \mu_\lambda(x_i) \,
                  \cauchy_{\lambda}^{\underline{\h_i-1}}(x_i-1),
              \end{equation}
              where
              \begin{equation}
                  \label{eq:g-is-a-magic-product}
                  \func(x_1, \dots, x_\ell) = \prod_{1 \leq i < j \leq \ell}
                  \frac{(x_i-x_j)(x_i-x_j-\h_i+\h_j)}{(x_i-x_j+\h_j)(x_i-x_j-\h_i)}.
              \end{equation}

    \end{enumerate}
    \begin{figure}
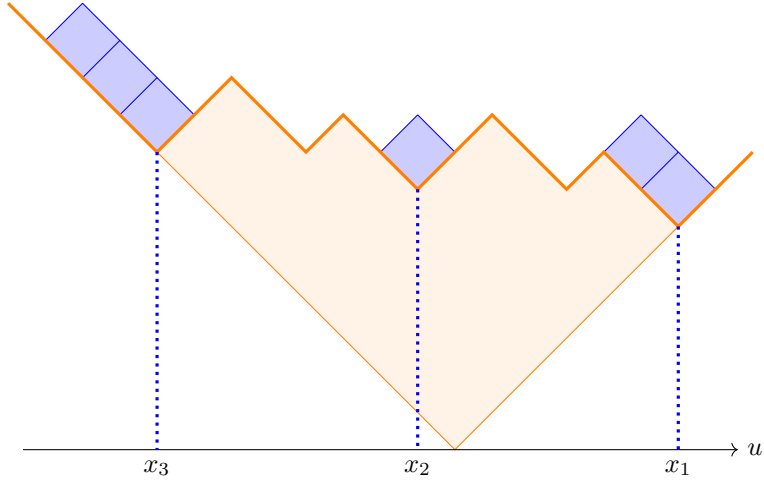

        \subfile{figures/FIGURE-Pieri.tex}
        \caption{Example of the growth of a Young
            diagram considered in \cref{lem:string}. The light orange area is the original
            Young diagram $\lambda$, and the dark blue boxes are added in successive steps,
            from right to left. The quantities from \cref{lem:string}\ref{item:B} are as
            follows: $a_1 = 2$, $a_2 = 1$, $a_3 = 3$.}
        \label{fig:pieri}
    \end{figure}
\end{lemma}
Note that assumption \ref{item:annoying-condition} ensures that we do not
evaluate the Cauchy transform $\cauchy_\lambda$ or the function $\func$ at a
singularity on the right-hand side of \eqref{eq:big-product}.

\begin{proof}

    \emph{Proof of part \ref{item:A}.}
    Our strategy is to construct a coupling on a single probability space containing both:
    \begin{itemize}
        \item a uniformly random Poissonized tableau with shape $\lambda$, and
        \item a Plancherel growth process starting from $\lambda$.
    \end{itemize}
    An additional minor challenge is that our model requires some conditioning.

    \medskip

    Let $w_1,\dots,w_{\nklatki+k}$ be a sequence of \iid random variables with
    the uniform distribution $U(0,1)$ and let $\yd{i}= \RSK(w_1,\dots,w_i)$;
    then \[ \yd{0}\nearrow \cdots \nearrow \yd{\nklatki+k}\] is the Plancherel
    growth process. For $i\in\{1,\dots,k\}$ we denote by
    \[{\mathcal{U}_i=u(\yd{\nklatki+i}\setminus\yd{\nklatki+i-1})}\]
    the
    $u$-coordinate of the place where the growth occurs.

    Clearly, the probability distribution of the Plancherel growth process
    $(\diagram_{\nklatki},\dots,\diagram_{\nklatki+k})$ starting at
    $\diagram_{\nklatki}=\lambda$ coincides with the \emph{conditional}
    probability distribution of its counterpart $\big( \yd{\nklatki},\dots,
        \yd{\nklatki+k} \big)$, under the condition that
    \mbox{$\yd{\nklatki}=\lambda$}. By \cref{lem:conditioning}, the probability
    distribution of the random Poissonized tableau $T$ from the statement of the
    lemma coincides with the \emph{conditional} probability distribution of the
    insertion tableau $\mathcal{T}:=P(w_1,\dots,w_{\nklatki})$, under the
    condition $\yd{\nklatki}=\lambda$.

    These observations imply that it is enough to prove equality between the
    conditional expectations
    \begin{equation}
        \label{eq:claim}
        \E \giventhat*{
            \big( F_{\mathcal{T}}(u_0) \big)^k }{ \sigma\big(\yd{\nklatki} \big) } = k!\
        \E\giventhat*{  \indicator \big\{ u_0\geq
            \mathcal{U}_1>\cdots>\mathcal{U}_k\big\} }{ \sigma\big( \yd{\nklatki} \big)
        },
    \end{equation}
    where $\sigma \big(\yd{\nklatki} \big)$ denotes the
    $\sigma$-algebra generated by the random Young diagram~$\yd{\nklatki}$, and
    $\indicator\{A\}$ denotes the indicator random variable which takes the
    value $1$ if the condition $A$ holds true, and $0$ otherwise.

    \medskip

    For a moment let us fix the values in the prefix $x_1,\dots,x_{\nklatki}$;
    the conditional probability
    \begin{align} \label{eq:magic-1}
        \Pro
                  & \giventhat*{ F_{\mathcal{T}}(u_0) > w_{\nklatki+1} > \cdots >
        w_{\nklatki+k} }{ \sigma(x_1,\dots,x_{\nklatki}) }                        \\ \nonumber & = \E
           \giventhat*{  \indicator \left\{ F_{\mathcal{T}}(u_0) > w_{\nklatki+1} >
        \cdots > w_{\nklatki+k} \right\} }{ \sigma(x_1,\dots,x_{\nklatki})  }     \\
        \nonumber & = \operatorname{vol} \left\{
        (x_{\nklatki+1},\dots,x_{\nklatki+k}) \in [0,1]^k \colon
        F_{\mathcal{T}}(u_0) > w_{\nklatki+1} > \cdots > w_{\nklatki+k} \right\}  \\
        \nonumber & = \frac{1}{k!} \left[ F_{\mathcal{T}}(u_0) \right]^k
    \end{align}
    is then directly related to the value of the random variable
    $F_{\mathcal{T}}(u_0)$.

    The event which appears on the left hand side of \eqref{eq:magic-1} can be
    alternatively reformulated in the language of the Young diagrams
    $\yd{\nklatki}\nearrow \cdots \nearrow \yd{\nklatki+k}$ as follows:
    \begin{equation}
        \label{eq:insertion-monotone}
        \left\{ F_T(u_0) >
        w_{\nklatki+1} > \cdots > w_{\nklatki+k}  \right\} = \left\{   u_0\geq
        \mathcal{U}_1>\cdots>\mathcal{U}_k \right\}.
    \end{equation}
    Indeed, the
    equivalence \[ F_T(u_0) > w_{\nklatki+1}  \iff u_0\geq \mathcal{U}_1 \] is a
    consequence of the definition of $F_T(u_0)$ while each of the equivalences
    \[ w_{\nklatki+i} > w_{\nklatki+i+1} \iff \mathcal{U}_i > \mathcal{U}_{i+1}
    \] is the content of the Row Bumping Lemma \cite[page~9]{Fulton1997}. Thus,
    by taking the appropriate conditional expectation of both sides of
    \eqref{eq:magic-1}, the desired equality \eqref{eq:claim} follows
    immediately. %

    \medskip

    For part \ref{item:B} we refer to \cref{fig:pieri}.

    \medskip

    \emph{Proof of part \ref{item:C}.} We start with the case when the
    probability on the left-hand side of \eqref{eq:big-product} is non-zero. For
    an illustration see \cref{fig:pieri}. For integers $j\in\{1,\dots,\ell\}$
    and $m\in\{0,\dots,a_j\}$ we define the Young diagram $\lambda^{[j,m]}$ as
    the diagram $\lambda$ with additional boxes, the $u$-coordinates of which
    form the following multiset
    \begin{multline*}
        \underbrace{x_1,x_1-1, \dots,
            x_1-a_1+1}_{\text{$a_1$ elements}}, \dots, \underbrace{x_{j-1},x_{j-1}-1,
        \dots, x_{j-1}-a_{j-1}+1}_{\text{$a_{j-1}$ elements}}, \\
        \underbrace{x_{j},x_{j}-1, \dots, x_{j}-m+1}_{\text{$m$ elements}}.
    \end{multline*}
    Note that $\lambda^{[j,a_j]}=\lambda^{[j+1,0]}$. With this
    notation, the event $\mathbf{U} = \mathbf{x}^{\underline{\mathbf{a}}}$
    holds if and only if the sequence
    $(\diagram_{\nklatki},\dots, \diagram_{\nklatki+k})$ is equal to
    \begin{multline}
        \label{eq:path-of-PGP}
        \big( \underbrace{\lambda^{[1,0]},
            \dots \lambda^{[1,a_1]}}_{\text{$a_1+1$ elements}}, \underbrace{
            \lambda^{[2,1]}, \dots, \lambda^{[2,a_2]}}_{\text{$a_2$ elements}}, \dots,
        \underbrace{ \lambda^{[\ell,1]}, \dots,
            \lambda^{[\ell,a_\ell]}}_{\text{$a_\ell $ elements}}  \big)= \\[2ex] \big(
        \underbrace{\lambda^{[1,0]}, \dots \lambda^{[1,a_1-1]}}_{\text{$a_1$
                elements}}, %
        \dots, \underbrace{ \lambda^{[\ell-1,0]}, \dots,
            \lambda^{[\ell-1,a_{\ell-1}-1]}}_{\text{$a_{\ell-1}$ elements}},
        \underbrace{ \lambda^{[\ell,0]}, \dots,
            \lambda^{[\ell,a_\ell]}}_{\text{$a_\ell+1 $ elements}} \big).
    \end{multline}
    It follows that we must compute the probability of a Plancherel growth
    process, starting from $\lambda=\lambda^{[1,0]}$, traversing the sequence of
    diagrams \eqref{eq:path-of-PGP} in its first $k$ steps. Given the Markovian
    nature of this process, we calculate the probability of each transition
    independently and then take their product. We shall consider two categories of
    transitions separately:
    \begin{itemize}
        \item Those in which a new box is added to one of the concave corners of the
              original diagram $\lambda$, and \item The remaining transitions.
    \end{itemize}

    \medskip

    \emph{Transition from $\lambda^{[j,0]}$ to $\lambda^{[j,1]}$:} The diagram
    $\lambda^{[j,0]}$ is derived from $\lambda$ by the addition of $j-1$
    rectangles (see \cref{fig:pieri}). Each rectangle, indexed by $i \in
        \{1,\dots,j-1\}$, is defined by four vertices:
    \begin{itemize}
        \item Bottom and top vertices:
              \begin{itemize}
                  \item $u$-coordinates: $x_i$ and $x_i-a_i+1$, respectively.
                  \item Function: Either remove a concave corner of $\lambda$ or create a new convex corner.
                  \item Effect on Cauchy transform: either remove an existing pole or add an additional zero;  
                  algebraically this corresponds to multiplying by a factor of the form
                        \begin{equation}
                            \label{eq:monomial}
                            \big( z- (\text{$u$-coordinate of the vertex}) \big).
                        \end{equation}
              \end{itemize}
        \item Right and left vertices:
              \begin{itemize}
                  \item $u$-coordinates: $x_i+1$ and $x_i-a_j$, respectively.
                  \item Function: Either remove a convex corner of $\lambda$ or create a new concave corner.
                  \item Effect on Cauchy transform: either remove an existing zero or add an additional pole; 
                  algebraically this corresponds to dividing  by a factor of the form \eqref{eq:monomial}
              \end{itemize}

    \end{itemize}
    It follows that \[ \cauchy_{\lambda^{[j,0]}}(z) = \cauchy_{\lambda}(z)
        \prod_{i\in\{1,\dots,j-1\}}
        \frac{(z-x_i)(z-x_i+a_i-1)}{(z-x_i-1)(z-x_i+a_i)}. \]
    Thus the
    transition probability from the diagram $\lambda^{[j,0]}$ to
    $\lambda^{[j,1]}$ is equal to the residue
    \begin{multline}
        \label{eq:trans-1}
        \operatorname{Res}_{x_j} \cauchy_{\lambda^{[j,0]}}(z)
        = \left( \operatorname{Res}_{x_j} \cauchy_{\lambda}(z) \right)
        \prod_{i\in\{1,\dots,j-1\}}
        \frac{(x_j-x_i)(x_j-x_i+a_i-1)}{(x_j-x_i-1)(x_j-x_i+a_i)} = \\
        \mu_{\lambda}(x_j) \prod_{i\in\{1,\dots,j-1\}}
        \frac{(x_j-x_i)(x_j-x_i+a_i-1)}{(x_j-x_i-1)(x_j-x_i+a_i)}.
    \end{multline}

    \medskip

    \emph{Transition from $\lambda^{[j,m]}$ to $\lambda^{[j,m+1]}$ for $m>0$.}
    The diagram $\lambda^{[j,m]}$ can be derived from $\lambda$ by adding $j$ rectangles.
    Consequently, a reasoning similar to the one previously discussed implies that:
    \[ \cauchy_{\lambda^{[j,m]}}(z) = \cauchy_{\lambda}(z)
        \left( \prod_{i\in\{1,\dots,j-1\}}
        \frac{(z-x_i)(z-x_i+a_i-1)}{(z-x_i-1)(z-x_i+a_i)} \right)
        \frac{(z-x_j)(z-x_j+m-1)}{(z-x_j-1)(z-x_j+m)}. \]
    It follows that the
    transition probability from the diagram $\lambda^{[j,m]}$ to
    $\lambda^{[j,m+1]}$ is equal to the residue
    \begin{multline}
        \label{eq:trans-2}
        \operatorname{Res}_{x_j-m} \cauchy_{\lambda^{[j,0]}}(z) =
        \cauchy_{\lambda}(x_j-m ) \times \\ \prod_{i\in\{1,\dots,j-1\}}
        \frac{(x_j-m -x_i)(x_j-m -x_i+a_i-1)}{(x_j-m -x_i-1)(x_j-m -x_i+a_i)} \cdot
        \frac{(-1) m }{m+1}.
    \end{multline}

    \medskip

    \emph{The product.} We consider the probability \eqref{eq:trans-1}
    multiplied with the product of \eqref{eq:trans-2} over all choices of
    $m\in\{1,\dots,a_j-1\}$. Due to the telescopic cancellations
    this whole
    product is equal to \[ \mu_{\lambda}(x_j)\
        \cauchy^{\underline{a_j-1}}_{\lambda} (-1)^{a_j-1}  \frac{1}{a_j}
        \prod_{i\in\{1,\dots,j-1\} } \frac{x_j-x_i}{x_j-x_i-a_j}\cdot
        \frac{x_j-x_i+a_i-a_j}{x_j-x_i+a_i}. \] By taking the product over all
    choices of $j\in\{1,\dots,\ell\}$ we recover the right-hand side of
    \eqref{eq:big-product}, as required.

    \bigskip

    \emph{We consider now the case when the probability on the left-hand side of
        \eqref{eq:big-product} is equal to zero.} This means that at least one of
    the diagrams in the sequence \eqref{eq:path-of-PGP} is not well-defined. Let
    $\lambda^{[j,m]}$ with $j\in\{1,\dots,\ell\}$ and $m\in\{1,\dots,a_j\}$ be
    the first entry of this sequence which is not well-defined. This may happen
    only if $x_j-m+1$ is the $u$-coordinate of a convex corner of $\lambda$
    hence $\cauchy_{\lambda}(x_j-m+1)=0$ and $m\geq 2$. As a consequence, one of
    the factors on the right-hand side of \eqref{eq:big-product} is equal to
    zero, as required.
\end{proof}

\section{Decomposition into simple fractions}
\label{sec:simple-fractions}

In this section we will decompose the product $\func$ defined in
\eqref{eq:g-is-a-magic-product} into a sum of simple fractions.

\medskip

A \emph{spine graph} with $\en\geq 1$ vertices is defined as a directed path
graph $\gS$ such that the set of its edges consists of $\en-1$ elements and is
of the form \[ \sE_{\gS}= \{(v_1, v_2), \dots, (v_{\en-1}, v_{\en})\}.\] Note
that the vertices $v_1, \dots, v_{\en}$ are all different; otherwise, the graph
would not be connected. We denote the set of all spine graphs with the vertex
set $\sV=\{1, \dots, \en\}$ by $\sS_{\en}$; obviously $|\sS_{\en}|=\en!$. An
example of a spine graph is shown in \cref{fig:exspine}.

\begin{figure}
    \subfile{figures/FIGURE-spine.tex} \caption{An example of a spine
        graph with $6$ vertices.} \label{fig:exspine}

    \bigskip \bigskip \bigskip

    \subfile{figures/FIGURE-multi-spine.tex} \caption{An example of a
        multi-spine graph with $9$ vertices and $4$ connected components.}
    \label{fig:exmulti-spine}
\end{figure}

\emph{A multi-spine graph} is defined as any directed graph such that each
component is a spine graph. In other words, a multi-spine graph is a forest of
spine graphs. We denote the set of all multi-spine graphs with the vertex set
$\sV = \{ 1, \dots, \en \}$ by $\sMS_{\en}$. An example of a multi-spine graph
is shown in \cref{fig:exmulti-spine}.

\begin{lemma} \label{lem:lemprod}

    Let $(\h_1, \dots \h_{\en})$ be a sequence of numbers which has the property
    that the sum of the entries of any non-empty subsequence is non-zero (this
    condition holds, for example, if $\h_1,\dots,\h_{\en}>0$ are all positive).

    Then the element $\func\in\R(x_1,\dots,x_{\en})$ of the field of rational
    functions defined in \eqref{eq:g-is-a-magic-product} can be written as the
    sum %
    \begin{equation}
        \label{eq:lemprod} 
        \func(x_1,\dots,x_\te)
        =  \sum\limits_{\gMS \in \sMS_{\en}} \frac{ \beta_{\gMS} }
        {\prod\limits_{(i,j)\in \sE_{\gMS} }\left(x_j-x_i+\h_i\right)}.
    \end{equation}
    Above, for any graph $\gMS \in  \sMS_{\en}$, the constant
    $\beta_{\gMS}$ is defined as
    \begin{equation}
        \label{eq:beta}
        \beta_{\gMS} =
        (-1)^{|V_{\gMS}|}\  \frac{ \prod\limits_{j=1}^{\en} \h_j }{
            \prod\limits_{\gS'}\left[ (-1) \cdot  \sum\limits_{i \in V_{\gS'}} \h_i
                \right] } ,
    \end{equation}
    where the product over $\gS'$ runs over all connected components of the
    graph~$\gMS$. 
\end{lemma} 
\begin{proof}
    To simplify the notation, we put \[ z_j=x_j-\h_j \] for each index
    $j \in \{ 1, \dots, \en\}$. %

    Let \[\macierz=\left[\frac{1}{x_i-z_j}\right]_{1 \leq i, j \leq \en }\] be
    the Cauchy matrix \cite{Schechter1959}. Its determinant, called the
    \emph{Cauchy determinant}, is given by the following product formula
    \cite{Schechter1959} $$\det \macierz =\frac{\prod\limits_{1 \leq i < j \leq
                \en } (x_i-x_j)(z_j-z_i)}{\prod\limits_{1 \leq i, j \leq \en} (x_j-z_i)}.$$
    The denominator of $\func$ differs from its counterpart in the Cauchy
    determinant only by the missing diagonal factors $x_j-z_i$ for $i=j$. Thus
    \[ \func = \left( \prod\limits_{j=1}^{\en} \left(x_j-z_j\right) \right) \det
        \macierz = \left( \prod\limits_{j=1}^{\en} \h_j \right) \det \macierz  . \]

    Using the definition of the determinant we express $\func$ as a sum over
    permutations %
    \[ \func =\left( \prod\limits_{j=1}^{\en} \h_j \right) \sum\limits_{\sigma
            \in \Sy{\nmoment}} \frac{(-1)^{\en-c(\sigma)}}{ \prod\limits_{i=1}^{\en}
        (x_{\sigma(i)}-z_i) }, \] where $c(\sigma)$ denotes the number of cycles of
    the permutation $\sigma$. We can treat each permutation $\sigma \in \Sy{n}$
    as a directed weighted graph with the vertex set $\sV_{\sigma}=\{1, \dots,
        \en \}$ and with the edge set
    \[ \sE_{\sigma}=\left\{ \left(1,
        \sigma(1)\right), \dots, \left(\en, \sigma(\en)\right)\right\}.\]
    We define
    the weight of an edge $e=\left(i, \sigma(i)\right)$ as $\weight(e) = a_{i}$.
    Consequently, we can express $\func$ as:
    \begin{equation}
        \label{eq:wzorek}
        \func = \left(
        \prod\limits_{j=1}^{\en} \h_j \right) \sum\limits_{\sigma \in \Sy{\nmoment}}
        (-1)^{\en-c(\sigma)} f_{\sigma} = (-1)^{ \en } \left(
        \prod\limits_{j=1}^{\en} \h_j \right)  \sum\limits_{\sigma \in
            \Sy{\nmoment}} \prod\limits_{\sigma'} \left( -f_{\sigma'} \right),
    \end{equation}
    where $\sigma'$ ranges over the connected components of the directed graph
    represented by $\sigma$. Each such connected component corresponds to a cycle
    in the permutation~$\sigma$. It is worth noting that $f_{\sigma}$ was
    previously defined in \cref{sec:rational}.

    Let $\sigma'$ be a connected component of the directed graph $\sigma$. Using the identity
    \[
        \sum\limits_{i \in \sV_{\sigma'}} \h_{i} = \sum\limits_{j \in \sV_{\sigma'}} x_{j} - \sum\limits_{i \in \sV_{\sigma'}} z_{i} =
        \sum\limits_{(i,j) \in \sE_{\sigma'}} \left( x_{j} - z_{i} \right)
    \]
    we obtain
    \begin{equation}
        \label{eq:cancel}
        f_{\sigma'} \sum\limits_{i \in \sV_{\sigma'}} \h_{i} = f_{\sigma'} \sum\limits_{(i,j) \in \sE_{\sigma'}} (x_j - z_i) = \sum\limits_{\gS'} f_{\gS'},
    \end{equation}
    where $\gS'$ runs over all spine graphs obtained from the cycle $\sigma'$ by removing exactly one edge.

    Equation \eqref{eq:cancel} can be written as
    \[
        f_{\sigma'} = \frac{1}{\sum\limits_{i \in \sV_{\sigma'}} \h_{i}} \sum\limits_{\gS'} f_{\gS'},
    \]
    we apply this identity to each cycle $\sigma'$ of the permutation $\sigma \in
        \Sy{\nmoment}$ on the right-hand side of \cref{eq:wzorek}. Note that the above
    equality holds true also in the special case when the cycle $\sigma'$ is a
    fix-point; in this case, we remove the loop from the directed graph $\sigma'$,
    and the unique resulting graph $\gS'$ has one isolated vertex and no edges.

    If we remove one edge from each cycle of every permutation in all possible
    ways, we obtain each multi-spine graph on the vertex set $\{1, \dots, \en\}$
    exactly once. In this way, we proved that
    \[
        \func = (-1)^{\en} \left( \prod\limits_{j=1}^{\en} \h_j \right) \sum\limits_{\gMS \in \sMS_{\en}} \prod\limits_{\gS'} \frac{- f_{\gS'}}{\sum\limits_{i \in \sV_{\gS'}} a_i},
    \]
    where $\gS'$ runs over the connected components of the graph $\gMS$, as required.
\end{proof}

\section{The moments of the threshold}
\label{rozdzial6}

\subsection{The first formula for the moments} A \emph{composition} of a natural
number $\nmoment$ is an expression of $\nmoment$ as an ordered sum of positive
integers $\nmoment=a_1+\cdots+a_{\te}$. The set of all compositions of
$\nmoment$ will be denoted by $\com_{\nmoment}$. For a given composition
$\mathbf{a}=(a_1,\dots,a_{\te}) \in \com_{\nmoment}$ we denote the number of its
parts by $\te=\te(\mathbf{a})$.

Using \cref{lem:string} we obtain
\begin{multline}
    \label{eq:moment-formula}
    m_{\nmoment}\left(F_T\left(u_0\right)\right) = \nmoment ! \sum_{\mathbf{\h} \in
        \comp_{\nmoment}} \sum_{\mathbf{x}} \Pro \big[ \mathbf{U} =
    \mathbf{x}^{\underline{\mathbf{a}}}\big] \\ = \nmoment ! \sum_{\mathbf{\h} \in
        \comp_{\nmoment}} \sum_{\mathbf{x}} \func(x_1,\dots,x_{\te}) \prod_{i=1}^{\te}
    \frac{ (-1)^{\h_i-1}}{\h_i} \ \mu_\lambda(x_i)\
    \cauchy^{\underline{\h_i-1}}_{\lambda}(x_i-1),
\end{multline}
where in each
expression the second sum runs over $\mathbf{x}=(x_1,\dots,x_{\te})\in\Deco$
such that %
\begin{equation}
    \label{eq:decreasing}
    u_0 \geq x_{1} > x_{2} > \cdots > x_\te
\end{equation}
\emph{and} such that the condition \ref{item:annoying-condition}
from \cref{lem:string} is satisfied.

\medskip

The condition \ref{item:annoying-condition} proves to be rather unwieldy in
practical applications. To address this, our strategy is to derive an analogue
of formula \eqref{eq:moment-formula} that involves summation over \emph{all}
tuples $x_1,\dots,x_{\te}\in\Deco$ satisfying \eqref{eq:decreasing}, effectively
eliminating the need for condition \ref{item:annoying-condition}. However,
without this additional constraint, there is a risk that one of the factors in
the falling product $\cauchy^{\underline{\h_i-1}}_{\lambda}(x_i-1)$ might be
evaluated at a singularity, potentially leading to division by zero in the
right-hand side of \eqref{eq:moment-formula}.

To circumvent this issue, we will expand our focus beyond Young diagrams to
consider a more general class of objects known as \emph{interlacing sequences}.
This approach allows us to easily avoid such divisions by zero. The formulas for
the Young diagram $\lambda$ can then be obtained through an appropriate limiting
process.

\subsection{Interlacing sequences}

The following notations are based on the work of Kerov \cite{Kerov1993}. We say
that
\begin{equation}
    \label{eq:this-is-an-interlacing}
    \Lambda=(\kerx_0,\dots,\kerx_\kerell; \; \kery_1,\dots,\kery_{\kerell})
\end{equation}
is an \emph{interlacing sequence}  if its entries are real
numbers such that
\[ \kerx_0 < \kery_1 < \kerx_1 < \cdots <  \kery_\kerell < \kerx_\kerell. \]
Following \cref{fig:cauchy} and \cref{sec:transition-measure}, each Young
diagram can be regarded as an interlacing sequence. Conversely, each
interlacing sequence can be visualized as a zig-zag curve similar to the one in
\cref{fig:cauchy}. Therefore, we will refer to the entries of the sequence
$\kerx_0, \dots, \kerx_\kerell$ as \emph{concave corners} and the entries of the
sequence $\kery_1, \dots, \kery_\kerell$ as \emph{convex corners}.

The Cauchy transform $\cauchy_\Lambda$ and the transition measure
$\mu_{\Lambda}$ of an interlacing sequence $\Lambda$ is defined in an analogous
way as their counterparts for Young diagrams in \cref{sec:transition-measure}.

\subsection{Moments for interlacing sequences} \label{sec:moments}

Let an interlacing sequence $\Lambda$ be fixed. We assume that the set of
concave corners is \emph{generic}, i.e., if $i\neq j$ then $\kerx_i-\kerx_j$ is
\emph{not} an integer. For the set of decoration values $\Deco:=\{ \kerx_0,
    \dots, \kerx_\kerell \}$ we take the concave corners. Let $u_0$ be a fixed real
number. We define the $\nmoment$-th \emph{moment} for the interlacing sequence
$\Lambda$ as
\begin{equation}
    \label{eq:mominterseq}
    \Moment_{\nmoment}=
    \Moment_{\nmoment} (\Lambda, u_0) =  \nmoment ! \sum_{\mathbf{\h} \in
        \mathcal{C}_n} \sum_{\mathbf{x}} \func(x_1,\dots,x_{\te}) \prod_{i=1}^{\te}
    \frac{ (-1)^{\h_i-1}}{\h_i} \ \mu_\Lambda(x_i)\
    \cauchy^{\underline{\h_i-1}}_{\Lambda}(x_i-1),
\end{equation}
where the sum over $\mathbf{x}$ runs over $x_1,\dots,x_{\te}\in\Deco$ such that
\eqref{eq:decreasing} holds true, and $\te=\te(\mathbf{a})$ denotes the length
of the composition $\mathbf{a}$ as before. The assumption of generic concave
corners ensures that the right-hand side is well-defined. One might question
whether $M_{\nmoment} (\Lambda, u_0)$ has a probabilistic interpretation as a
moment of some natural random variable associated with the interlacing sequence
$\Lambda$. We conjecture that this is not the case. Instead, we will employ
$M_{\nmoment}$ purely as an auxiliary tool for investigating the moments of the
random variable $F_T(u_0)$, as discussed below.

The right-hand side of \eqref{eq:mominterseq} is very similar to its counterpart
\eqref{eq:moment-formula}; the only difference is that the second sum on the
right-hand side of \eqref{eq:moment-formula} runs over certain
sequences~$\mathbf{x}$ which \emph{additionally} fulfill the condition
\ref{item:annoying-condition} from \cref{lem:string}\ref{item:C}.

Let us fix an integer $s\in\{0,\dots,\kerell+1\}$ and consider the set
$W_{s,\kerell}$ of\linebreak[4] interlacing sequences $\Lambda$ of the form
\eqref{eq:this-is-an-interlacing} with the property that
${\kerx_0,\kerx_1,\dots,\kerx_{s-1} \leq u_0 }$ are all small and
$\kerx_s,\dots,\kerx_{\kerell} > u_0$ are all big; in other words $s$ is the
cardinality of small entries of the set $\Deco$. Thanks to the aforementioned
removal of the condition \ref{item:annoying-condition}, the restriction of the
function $\Lambda\mapsto\Moment_{\nmoment}(\Lambda,u_0)$ to the set
$W_{s,\kerell}$ is a rational function in the variables
$\kerx_0,\dots,\kerx_\kerell,\kery_1,\dots,\kery_\kerell$. Our general strategy
is to investigate this rational function $M_{\nmoment}$.

The omission of condition \ref{item:annoying-condition} results in the rational
function $\Moment_{\nmoment}$ being singular for certain non-generic interlacing
sequences. Specifically, it is unclear how to evaluate
$\Moment_{\nmoment}(\Lambda, u_0)$ when the interlacing sequence $\Lambda$
corresponds to a Young diagram~$\lambda$, which is inherently non-generic.
However, \cref{cor:regularize} demonstrates a special method for taking the
limit of $\Moment_{\nmoment}$ at the singularity, thereby connecting it to our
primary focus: the moment $m_{\nmoment}\left(F_T\left(u_0\right)\right)$.

Furthermore, the proof of \cref{thm:wzormikolaja} will later reveal that this
singularity is removable. Consequently, an analogue of \cref{cor:regularize}
holds true for \emph{any} method of taking the limit $\Lambda \to \lambda$.

\subsection{Regularization}

Let a Young diagram $\lambda$ be fixed and let $\Lambda$ be the corresponding
interlacing sequence. For $\epsilon>0$ we define the interlacing sequence
\[
    \Lambda^\epsilon=(\kerx^\epsilon_0,\dots,\kerx^\epsilon_{\kerell}; \;
    \kery^\epsilon_1,\dots,\kery^\epsilon_{\kerell})\]
given by
\[ \kerx^\epsilon_j
    = \kerx_j + j \epsilon, \qquad \kery^\epsilon_j = \kery_j + j \epsilon. \]
Note
that if $\epsilon$ is small enough, the set of concave corners of
$\Lambda^\epsilon$ is generic so that $ \Moment_{\nmoment} (\Lambda^\epsilon,
    u_0)$ is well-defined.

The distance
\begin{equation}
    \label{eq:distance}
    \kerx^\epsilon_j-
    \kery^\epsilon_j = \kerx_j- \kery_j
\end{equation}
between any convex corner
$\kery^\epsilon_j$ and the next concave corner to the right $\kerx^\epsilon_j$
does not depend on the value of $\epsilon$, and is a positive integer which has
a natural interpretation for the original Young diagram $\lambda$,
cf.~\cref{fig:cauchy}.

\begin{lemma} \label{cor:regularize}

    We suppose that $u_0$ is \emph{not} an integer number. With the above
    notations, the moment $m_{\nmoment}$ is equal to the limit of the moment
    $\Moment_{\nmoment}$, when $\epsilon$ tends to zero: \[
        m_{\nmoment}\left(F_T\left(u_0\right)\right) =  \lim_{\epsilon\to 0}
        \Moment_{\nmoment} (\Lambda^\epsilon, u_0). \] \end{lemma}
\begin{proof}
    Let
    $s\in\{0,\dots,\kerell+1\}$ be the cardinality of the small concave corners
    of $\lambda$; with the notations of \cref{sec:moments} this means that
    $\Lambda^\epsilon\in W_{s,\kerell}$ if $|\epsilon|$ is small enough. By
    writing $x_i=\kerx_{r_i}^\epsilon$ we may write  \eqref{eq:mominterseq} as
    \begin{equation}
        \label{eq:mominterseq=RR}
        \Moment_{\nmoment}(\Lambda^\epsilon,u_0) =  \nmoment ! \sum_{\mathbf{\h} \in
            \mathcal{C}_n} \sum_{s\geq r_1 > \cdots > r_\te \geq 1}
        \func(\kerx_{r_1}^\epsilon,\dots,\kerx_{r_\te}^\epsilon) \prod_{i=1}^{\te}
        \frac{ (-1)^{\h_i-1}}{\h_i} \ \mu_\Lambda(\kerx_{r_i}^\epsilon)\
        \cauchy^{\underline{\h_i-1}}_{\Lambda^\epsilon}(\kerx_{r_i}^\epsilon-1).
    \end{equation}
    Similarly \eqref{eq:moment-formula} can be written as
    \begin{equation}
        \label{eq:mominterseq=SS}
        m_{\nmoment}\left(F_T\left(u_0\right)\right) =  \nmoment ! \sum_{\mathbf{\h}
            \in \mathcal{C}_n} \sum \func(\kerx_{r_1},\dots,\kerx_{r_\te})
        \prod_{i=1}^{\te} \frac{ (-1)^{\h_i-1}}{\h_i} \ \mu_\Lambda(\kerx_{r_i})\
        \cauchy^{\underline{\h_i-1}}_{\Lambda}(\kerx_{r_i}-1);
    \end{equation}
    the
    consequence of the condition \ref{item:annoying-condition} from
    \cref{lem:string} is that the second sum runs over $s\geq r_1 > \cdots >
        r_\te \geq 1$ which additionally fulfill
    \begin{equation}
        \label{eq:XX}
        \kerx_{r_i} - \kery_{r_i} \geq a_i \qquad \text{for } i\in\{1,\dots,\te\} ;
    \end{equation}
    in the special case when $r_i=0$ and $\kery_0$ is not defined
    the above condition is fulfilled by convention.

    Let us consider a summand of \eqref{eq:mominterseq=RR} which corresponds to
    $\mathbf{\h} \in \mathcal{C}_n$ and a tuple $(r_1,\dots,r_\te)$ for which
    \eqref{eq:XX} is \emph{not} satisfied thus  $1\leq \kerx_{r_i} - \kery_{r_i}
        \leq a_i-1$ for some choice of the index $i$. One of the factors in
    $\cauchy^{\underline{\h_i-1}}_{\Lambda^\epsilon}(\kerx_{r_i}^\epsilon-1)$ is
    equal to \[ \cauchy_{\Lambda^\epsilon}\big( \kerx_{r_i}^\epsilon-
        (\kerx_{r_i}-\kery_{r_i} ) \big) =  \cauchy_{\Lambda^\epsilon}\big(
        \kerx_{r_i}^\epsilon- (\kerx_{r_i}^\epsilon-\kery_{r_i}^\epsilon \big) =
        \cauchy_{\Lambda^\epsilon}\big( \kery_{r_i}^\epsilon \big) =0 \] by the very
    definition of the Cauchy transform; as a consequence the whole corresponding
    summand of \eqref{eq:mominterseq=RR} vanishes as well.

    On the other hand, any summand in \eqref{eq:mominterseq=RR} for which
    \eqref{eq:XX} is satisfied is continuous at $\epsilon=0$ and clearly
    converges as $\epsilon\to 0$ to its counterpart in \eqref{eq:mominterseq=SS}
    which completes the proof.
\end{proof}

\subsection{Cumulants for interlacing sequences}
\label{sec:cumulants-interlacing}

For a given interlacing sequence $\Lambda$ and $u_0$ we consider the
corresponding sequence of moments $M_1,M_2,\dots$ with $
    M_\nmoment=M_\nmoment(\Lambda,u_0)$ given by \eqref{eq:mominterseq}. We revisit
\cref{sec:cumulants-and-moments} and consider the corresponding sequence of
formal cumulants $K_1,K_2,\dots$ with $K_k=K_k(\Lambda,u_0)$ given by the
expansion $$\log \sum_{\nmoment=0}^{\infty} \frac{M_{\nmoment}}{\nmoment !}
    t^{\nmoment}  = \sum_{\nmoment=1}^{\infty} K_{\nmoment}
    \frac{t^{\nmoment}}{\nmoment !}.$$

Since each cumulant $K_k$ can be expressed as a polynomial in the moments
$M_1,\dots,M_k$,  \cref{cor:regularize}  implies the following result.

\begin{lemma}
    \label{cor:regularize2}

    Suppose that $u_0$ is \emph{not} an integer number. With the above notations,
    the cumulants of the random variable $F_T\left(u_0\right)$ are given by
    \[
        \kappa_{\nmoment}\left(F_T\left(u_0\right)\right) =  \lim_{\epsilon\to 0}
        K_{\nmoment} (\Lambda^\epsilon, u_0). \]
\end{lemma}

\section{Proof of \cref{thm:wzormikolaja}} \label{sec:dowod-kumulant-konczy-sie}

The current section is devoted to the proof of \cref{thm:wzormikolaja}. The key
ideas are as follows.

In \crefrange{sec:multi-caterpillar-graphs}{sec:black-decreasing-decorations},
we find a formula \eqref{eq:kontynuuj} for the moment
$\Moment_{\nmoment}(\Lambda,u_0)$, which is expressed in terms of
\emph{multi-caterpillar graphs} (the set $\sMC_{\nmoment}^{\tagg}$) and their
\emph{black-decreasing decorations}. These multi-caterpillar graphs come in
three flavors, differing by the way their vertices are indexed, which affects
the way they are counted (see \cref{sec:three-systems}).

By the first application of the double counting technique in
\cref{sec:double-counting-1}, we transform this formula for
$\Moment_{\nmoment}(\Lambda,u_0)$ to another one, expressed in terms of the
multi-caterpillar graphs from the set $\sMC_{\nmoment}^{\markk}$, which have
better symmetry properties, and their \emph{black-injective decorations}. In
\cref{sec:no-injective}, we show that in the latter sum, the restriction to
black-injective decorations can be lifted, which simplifies matters.

In \cref{sec:first-cumulant-formula}, we use the fact that the cumulants
correspond to taking the contribution from only the connected graphs, and we
obtain an intermediate formula for the cumulants.

In
\crefrange{sec:cat-labeled}{sec:second-formula-cumulants}, we apply the double
counting technique for the second time and obtain a formula for the cumulant,
this time in terms of the caterpillar graphs from the set
$\sC^{\labeled}_{\nmoment}$ and their decorations.

In \cref{sec:sum-ncat}, we
transform this formula into a final form, with the sum over non-crossing
alternating trees and their decorations.

Finally, in \cref{sec:proof-limit}, we
use the link between the cumulants for Young diagrams and for interlacing
sequences to complete the proof.

\subsection{The graph expansion for the moments}

Using \cref{lem:lemprod} and the fact that for any integer $r\geq 1$
\begin{equation}
    \cauchy_\Lambda(x_i-r) = -\sum_{x_{i,r} }
    \frac{\mu_{\Lambda}(x_{i,r})}{x_{i,r}-x_i+r},
\end{equation}
we may rewrite the
formula   \eqref{eq:mominterseq} \label{eq:cytpow2}
as follows.

\begin{corollary}
    \label{cor:formula-for-moments}

    If the interlacing sequence $\Lambda$ is generic then the moment $M_k$ is given
    by
    \begin{multline}
        \label{eq:moment-composition}
        \Moment_{\nmoment}(\Lambda,u_0) = \nmoment !\ \sum_{\mathbf{\h} \in
            \com_{\nmoment}} \sum_{x_{\te} < \cdots < x_1 \leq u_0}  \prod_{i=1}^{\te}
        \frac{\mu_{\Lambda}(x_i)}{\h_i}  \times \\ \prod_{r=1}^{\h_i-1} \sum_{x_{i,r} }
        \frac{\mu_{\Lambda}(x_{i,r})}{x_{i,r}-x_i+r}  \times \sum\limits_{\gMS \in
            \sMS_{\en}} \frac{\beta_{\gMS}} {\prod\limits_{(i,j)\in
                \sE_{\gMS}}\left(x_j-x_i+\h_i\right)},
    \end{multline}
    where
    $\ell=\ell(\mathbf{\h})$ is the length of the composition $\mathbf{\h}$.
    Recall that the constant $\beta_{\gMS}$ was defined in \eqref{eq:beta}. The
    above sums run over $x_i , x_{i,r} \in \{\kerx_1,\dots,\kerx_{\kerell} \}$.
\end{corollary}

In the following we denote \[x_{i,0}:=x_i.\]

Now we will define multi-caterpillar graphs and with them we will simplify
\cref{cor:formula-for-moments}.

\subsection{Multi-caterpillar graphs}
\label{sec:multi-caterpillar-graphs}

By applying the distributive law to the right-hand side of
\eqref{eq:moment-composition} we obtain a sum of a lot of terms; to each of them
we shall associate a certain directed weighted graph $\gC$. Each term is a
product of:
\begin{itemize}
    \item the numerical factor \[ \nmoment!\ \beta_{\gMS}
              \prod\limits_{i=1}^{\te} \frac{1}{\h_i} \prod\limits_{r\in\{0,\dots, \h_i-1\}}
              \mu_{\Lambda}(x_{i,r})  \] for some multi-spine graph $\gMS$, and
    \item the
          reciprocal of the product of the polynomials of the form \[ (x_{i,r}-x_{i,0}+r)
              \quad \text{ or } \quad  (x_{j,0}-x_{i,0}+\h_j).\]
\end{itemize}
The latter
product of polynomials is in our focus.

\begin{figure}
    \subfile{figures/FIGURE-graph-on-diagram}

    \caption{The multi-caterpillar graph considered in \cref{example:caterpillar}.
        The composition $\mathbf{a} = (1,2,2,1,3)$ is visualized as a configuration of
        white boxes representing potential anti-Pieri growth. The vertices of the
        multi-caterpillar graph (black vertices $\bullet$ and red vertices $\otimes$)
        correspond to the boxes of the Young diagram where the Plancherel growth
        occurred. For clarity, the weights of the edges are not shown.}
    \label{fig:dwaraz}

    \medskip

    \subfile{figures/FIGURE-graph-without-diagram}

    \caption{The multi-caterpillar graph from \cref{fig:dwaraz} with the weights of
        the edges shown. }
    \label{fig:dwadwa}
\end{figure}

\begin{example}
    \label{example:caterpillar}
    The graph shown in \cref{fig:dwaraz}
    was obtained from the term
    \begin{multline*}
        \frac{1}{(x_{5,1}-x_5+1)(x_{5,2}-x_5+2)(x_{3,1}-x_3+1)(x_{2,1}-x_2+1)} \times \\
        \frac{1}{(x_5-x_1+1)(x_4-x_5+3)(x_3-x_4+1)}
    \end{multline*}
    which is one of the
    summands in \cref{cor:formula-for-moments} which corresponds to
    ${\mathbf{a}=(1,2,2,1,3).}$ \cref{fig:dwadwa} shows the same graph without the
    Young diagram.
\end{example}

\subsection{Multi-caterpillar graphs, the formal approach}
\label{sec:multi-caterpillar-formal}

More formally, a \emph{multi-caterpillar graph $G$ with tagged vertices} is a
directed, weighted graph containing black and red vertices that satisfies the
following properties:
\begin{itemize}
    \item There exists a tuple of integers $\h_1, \ldots, \h_\ell \geq 1$, where $\ell \geq 1$.

    \item The subgraph induced by the black vertices forms a multi-spine graph with $\ell$ vertices tagged
          \[(1,0), \ldots, (\ell,0).\]
          For each $j \in \{1, \ldots, \ell\}$, the black vertex tagged $(j,0)$ has at most one outgoing edge to
          another black vertex, with weight $\h_j$ if such an edge exists.

    \item Removing the edges between black vertices results in $\ell$ connected components. Each component consists of:
          \begin{itemize}
              \item a single black vertex tagged $(j,0)$ for some $j \in \{1, \ldots, \ell\},$
              \item and $\h_j - 1$ red vertices tagged $(j,1), \ldots, (j,\h_j - 1)$.
          \end{itemize}

    \item Within each component:
          \begin{itemize}
              \item There are no edges between red vertices.
              \item For each $k \in \{1, \ldots, \h_j - 1\}$, there is a directed edge from the black vertex $(j,0)$ to the red vertex $(j,k)$ with weight $k$.
          \end{itemize}
\end{itemize}

Let $\sMC_{\nmoment}^{\tagg}$ denote the set of all multi-caterpillar graphs with $n$ tagged vertices.

\subsection{Three systems of vertex nomenclature}
\label{sec:three-systems}

In the following discussion,
we will employ the technique of double counting twice. To facilitate this, we
introduce three distinct systems for naming the vertices in a directed weighted
graph with $\nmoment$ vertices:
\begin{itemize}
    \item \textbf{Tags}: Elements of the set $\N \times \N_0$. This tagging system
          was utilized in
          \cref{sec:multi-caterpillar-graphs,sec:multi-caterpillar-formal}.
          \protect\cref{fig:przyka} illustrates a graph with tagged vertices.

    \item \textbf{Labels}: Elements of the set $\{1,\dots,\nmoment\}$. We consider
          only labelings that satisfy the following property: For any directed edge
          $e=(v_1,v_2)$ connecting vertices $v_1,v_2 \in \{1,\dots,\nmoment\}$, the
          weight of the edge is equal to the difference of the vertex labels:
          \begin{equation}
              \label{eq:labeling-compatible}
              \weight(e) = \weight(v_1,v_2) = v_2 - v_1.
          \end{equation}
          This is consistent with \cref{eq:konwencja-wag}.

    \item \textbf{Marks}: Elements of an arbitrary fixed set containing $\nmoment$
          elements. To distinguish marks from labels, we may define the set of marks as:
          \begin{equation}
              \label{eq:underlined}
              \{ \underline{1}, \underline{2}, \dots, \underline{\nmoment} \}
          \end{equation}
          consisting of $\nmoment$ underlined integers. \cref{fig:przykb} demonstrates a
          graph with marked vertices.

\end{itemize}

\subsection{Black-decreasing decorations}
\label{sec:black-decreasing-decorations}

Let $\gMC$ be a multi-caterpillar graph with tagged vertices. A decoration
$\mathbf{x} \in D_{\gMC}$ is called \emph{black-decreasing} if, for any pair of
black vertices $(p,0)$ and $(q,0)$ with $p < q$, the corresponding values of the
decoration satisfy $x_{p,0} > x_{q,0}$. The set of all black-decreasing
decorations of a multi-caterpillar graph $\gMC$ will be denoted by
$D_{\gMC}^{>}$.

Using \cref{cor:formula-for-moments}, we can express the moment
$\Moment_{\nmoment}(\Lambda, u_0)$ as a sum over multi-caterpillar graphs. We
replace the double sum in \eqref{eq:moment-composition} over compositions and
multi-spine graphs with a sum over multi-caterpillar graphs $\gMC \in
    \sMC_{\nmoment}^{\tagg}$. Additionally, we replace the sum over the variables
$(x_i)$ and $(x_{i,r})$ with a sum over black-decreasing decorations. This gives
us:
\begin{equation}
    \label{eq:kontynuuj}
    \Moment_{\nmoment} =
    \nmoment ! \sum_{\gMC \in \sMC_{\nmoment}^{\tagg}} \sum_{\mathbf{x} \in D_{\gMC}^{>}} \alpha_{\gMC} f_{\gMC},
\end{equation}
where the constant $\alpha_{\gMC}$ is defined as:
\begin{equation}
    \label{eq:alphadef}
    \alpha_{\gMC} =
    (-1)^{|B_{\gMC}|} \left( \prod_{(i,j) \in V_{\gMC}} \mu_{\Lambda}(x_{i,j}) \right) \left( \prod_{\gC'} \frac{-1}{|V_{\gC'}|} \right),
\end{equation}
with $\gC'$ running over all connected components of the graph $\gMC$. Note that
$\alpha_{\gMC}$ also depends on the choice of the decoration $\mathbf{x}$; to
keep the notation lightweight, we will make this dependence implicit.

\subsection{The first double counting}
\label{sec:double-counting-1}

\begin{figure}
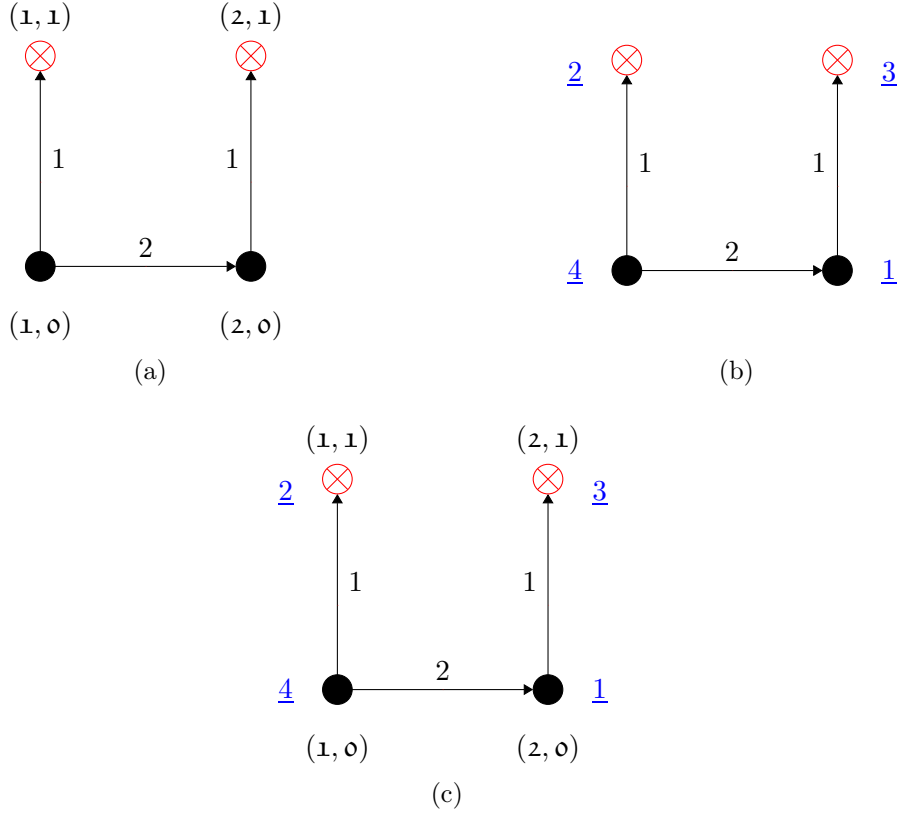

    \subfile{figures/FIGURE-tagged-and-marked.tex}
    \caption{Examples of multi-caterpillar graphs:
        \protect\subref{fig:przyka} A multi-caterpillar graph with tagged vertices.
        The tags, printed in black, are elements of the set $\left\{ \mathfrak{(1,0)},
            \mathfrak{(1,1)}, \mathfrak{(2,0)}, \mathfrak{(2,1)} \right\}$.
        \protect\subref{fig:przykb} A multi-caterpillar graph with marked vertices.
        The marks, printed in blue, are elements of the set $\{
            \underline{1},\underline{2},\underline{3},\underline{4}\}$.
        \protect\subref{fig:przykc} A multi-caterpillar graph with both tagged and
        marked vertices.
    }
    \label{fig:przyklady}
\end{figure}

Let $\sMC^{\markk}_{\nmoment}$ denote \emph{the set of multi-caterpillar graphs
    with $\nmoment$ marked vertices}, i.e., the set of weighted and directed graphs
$\gMC$ with vertex set $\{\underline{1}, \dots, \underline{\nmoment}\}$ such
that there exists a way to tag the vertices so that $\gMC$ becomes \emph{a
    multi-caterpillar graph with $\nmoment$ tagged vertices} as defined in
\cref{sec:multi-caterpillar-formal}. Let $\sMC_{\nmoment}^{\both}$ denote the
set of multi-caterpillar graphs with $\nmoment$ vertices that are simultaneously
tagged and marked. Examples of such graphs are shown in \cref{fig:przyklady}.

For any graph $\gMC \in \sMC_{\nmoment}^{\tagg}$, there are $\nmoment!$ ways to
mark its $\nmoment$ vertices by the elements of \eqref{eq:underlined}.

\medskip

Let $\gbi$ be a graph. Its decoration $\mathbf{x} \in D_{\gbi}$ is called
\emph{black-injective} if $x_i \neq x_j$ for all pairs of black vertices $i,j
    \in \sB_{\gbi}$ such that $i \neq j$. We denote the set of all black-injective
decorations of $\gbi$ by $D_{\gbi}^{\neq}$, and the set of non-black-injective
decorations of $\gbi$ by $D_{\gbi}^{=} = D_{\gbi} \setminus D_{\gbi}^{\neq}$.

Moreover, for each black-injective decoration $\mathbf{x}$ of $\gMC \in
    \sMC^{\markk}_{\nmoment}$, we can tag the vertices of $\gMC$ canonically as
follows:
\begin{itemize}
    \item Tag the black vertices by $(1, 0), (2, 0), \dots$ in the opposite order
          to that given by the decoration $\mathbf{x}$.

    \item For each black vertex $b$ with tag $(j, 0)$, tag the white vertices
          connected to $b$ by $(j, 1), (j, 2), \dots$ according to the increasing order of
          their corresponding edge weights.
\end{itemize}
This process transforms $G$ into a caterpillar graph with tagged vertices, and
$\mathbf{x}$ becomes a decreasing decoration.

Using these observations and \cref{eq:kontynuuj}, we obtain:
\begin{equation}
    \label{eq:wzorposredni}
    \Moment_{\nmoment} = \nmoment ! \sum_{
        \gMC\in\sMC_{\nmoment}^{\tagg}} \sum_{\mathbf{x} \in D_{\gMC}^{>}}
    \alpha_{\gMC} f_{\gMC} %
    = \sum_{ \gMC\in\sMC_{\nmoment}^{\both}} \sum_{\mathbf{x} \in D_{\gMC}^{>}}
    \alpha_{\gMC} f_{\gMC} %
    = \sum_{ \gMC\in\sMC_{\nmoment}^{\markk}}  \sum_{\mathbf{x} \in D_{\gMC}^{\neq}}
    \alpha_{\gMC} f_{\gMC}. %
\end{equation}
The constant $\alpha_{\gMC}$ is defined in \eqref{eq:alphadef}.

\subsection{Removal of injectivity requirement}
\label{sec:no-injective}

\begin{proposition}
    \label{cor:next-formula-for-moments}
    The following double sum over
    multi-caterpillar graphs and their decorations remains the same when we
    restrict the sum to black-injective decorations, i.e., for each integer $\nmoment
        \geq 1$, we have:
    \begin{equation}
        \label{eq:injective-is-boring}
        \sum_{\gMC \in \sMC_{\nmoment}^{\markk}} \sum_{\mathbf{x} \in D_{\gMC}} \alpha_{\gMC} f_{\gMC} =
        \sum_{\gMC \in \sMC_{\nmoment}^{\markk}} \sum_{\mathbf{x} \in D_{\gMC}^{\neq}} \alpha_{\gMC} f_{\gMC}.
    \end{equation}
\end{proposition}

\begin{proof}
    Consider the difference between the left-hand side and the right-hand side of \eqref{eq:injective-is-boring}:
    \begin{equation}
        \label{eq:delta}
        \Delta = \sum_{\gMC \in \sMC_{\nmoment}^{\markk}} \sum_{\mathbf{x} \in D_{\gMC}^{=}} \alpha_{\gMC} f_{\gMC}.
    \end{equation}
    Our goal is to prove that $\Delta = 0$.

    Let $\nmoment$ be a fixed natural number, and let $\sB = \{ b_1, \dots, b_l \}
        \subseteq \{ \underline{1}, \dots, \underline{\nmoment} \}$ be a fixed set. Let
    $\sMS(\sB)$ denote the set of all multi-spine graphs $\gMS$ with vertex set
    $V_{\gMS} = \sB$. In particular, $\sMS(\{ \underline{1}, \dots,
        \underline{\nmoment} \}) = \sMS_{\nmoment}$. Let
    $\grafydek^{\sMC}_{\nmoment}(\sB)$ denote the set of all multi-caterpillar
    graphs $\gMCempty \in \sMC^{\markk}_{\nmoment}$ with $\nmoment$ marked vertices
    $\underline{1}, \dots, \underline{\nmoment}$ such that:
    \begin{itemize}
        \item The set of black vertices of $\gMCempty$ is given by $B_{\gMCempty} = \sB$.
        \item There is no edge in $\gMCempty$ connecting two black vertices.
    \end{itemize}

    \medskip

    Let $k$ be a fixed natural number. Every multi-caterpillar graph $\gMC \in
        \sMC^{\markk}_{\nmoment}$ can be uniquely decomposed into the union of two
    graphs: $\gMCempty \in \grafydek^{\sMC}_{k}(B_{\gMC})$ and $\gMS \in
        \sMS(B_{\gMC})$. In other words, the graph $\gMS$ consists of all black
    vertices of $\gMC$ and the edges between them, while $\gMCempty$ consists of
    all vertices of $\gMC$ and the remaining edges. Furthermore, for each vertex $v
        \in V_{\gMS}$, we define the number $a_v$ as the number of vertices in the
    connected component of $\gMCempty$ that contains $v$. Using the notation from
    \cref{lem:lemprod}, the constant $\beta_{\gMS}$ given by \eqref{eq:beta} is:
    \begin{equation}
        \label{eq:beta2}
        \beta_{\gMS} = (-1)^{|V_{\gMS}|} \frac{\prod_{v \in \sV_{\gMS}} \h_v}{\prod_{\gC'} \left[ (-1) |V_{\gC'}| \right]},
    \end{equation}
    where the product over $\gC'$ runs over all connected components of the graph
    $\gMC$. Additionally, for each edge $e = (i, j) \in \sE_{\gMS}$, we define its
    weight as $\weight(e) = a_i$.

    Now we define the constant $\gamma_{\gMCempty}$ such that:
    \[
        \alpha_{\gMC} = \beta_{\gMS} \; \gamma_{\gMCempty}.
    \]
    From \eqref{eq:alphadef} and \eqref{eq:beta2}, we obtain:
    \[
        \gamma_{\gMCempty} = \frac{\alpha_{\gMC}}{\beta_{\gMS}} = \frac{\prod_{v \in V_{\gMCempty}} \mu_{\Lambda}(x_{v})}{\prod_{v \in B_{\gMCempty}} a_v},
    \]
    which depends only on the graph $\gMCempty$ and the decoration $\mathbf{x}$.

    For any set $B \subseteq \{ \underline{1}, \dots, \underline{\nmoment} \}$, the
    union of each pair of graphs $\gMCempty \in \grafydek^{\sMC}_{k}(B)$ and $\gMS
        \in \sMS(B)$ as above is a multi-caterpillar graph with $k$ marked vertices.
    Therefore, we can replace the sum in \eqref{eq:delta} over all
    multi-caterpillar graphs with marked vertices by a triple sum over all possible
    sets of black vertices, over multi-caterpillar graphs, and over all multi-spine
    graphs. It follows that:
    \begin{align*}
        \Delta & = \sum_{B \subseteq \{ \underline{1}, \dots, \underline{\nmoment} \}} \sum_{\substack{\gMC \in \sMC^{\markk}_{\nmoment}                                                                                                                           \\ B_{\gMC} = B}} \sum_{\mathbf{x} \in D_{\gMC}^{=}} \alpha_{\gMC} f_{\gMC} \\
               & = \sum_{B \subseteq \{ \underline{1}, \dots, \underline{\nmoment} \}} \sum_{\gMCempty \in \grafydek^{\sMC}_{\nmoment}(B)} \sum_{\mathbf{x} \in D_{\gMCempty}^{=}} \gamma_{\gMCempty} f_{\gMCempty} \sum_{\gMS \in \sMS(B)} \beta_{\gMS} f_{\gMS}.
    \end{align*}

    Let $B = \{b_1, \dots, b_l\} \subseteq \{ \underline{1}, \dots,
        \underline{\nmoment} \}$ be a fixed set and $\mathbf{x}_B = (x_{b_1}, \dots,
        x_{b_l})$ be a fixed non-black-injective decoration of $B$. Using the formula
    \eqref{eq:lemprod}, we obtain that the internal sum is:
    \begin{equation*}
        \sum_{\gMS \in \sMS(B)} \beta_{\gMS} f_{\gMS}(x_{b_1}, \dots, x_{b_l}) = \func(x_{b_1}, \dots, x_{b_l}) = 0,
    \end{equation*}
    since at least one of the factors in the numerator of $\func$ is zero (see
    \eqref{eq:g-is-a-magic-product}). Thus, $\Delta = 0$ as required.
\end{proof}

\subsection{The first formula for the cumulants}
\label{sec:first-cumulant-formula}

We denote by $\sC^{\markk}_{\nmoment} \subset \sMC^{\markk}_{\nmoment}$ the set
of \emph{connected multi-caterpillar graphs with $\nmoment$ marked vertices}.
Its elements will be called \emph{caterpillar graphs with $\nmoment$ marked
    vertices}.

Using \cref{cor:next-formula-for-moments} we transform the formula
\eqref{eq:wzorposredni} to $$\Moment_{\nmoment}= \sum_{
        \gMC\in\sMC_{\nmoment}^{\markk}} \sum_{\mathbf{x} \in D_{\gMC}} \alpha_{\gMC}
    f_{\gMC}.$$ We can look separately at each connected component $\gC'$ of a
multi-caterpillar graph $\gMC$. The connected components correspond to the
blocks of a set-partition. Thus %
\begin{align} \label{eq:moment-cumulant2}
    \Moment_{\nmoment} & = \sum_{
        \gMC\in\sMC_{\nmoment}^{\markk}} \sum_{\mathbf{x} \in D_{\gMC}} \alpha_{\gMC}
    f_{\gMC}                                              \\ \nonumber & = \sum_{ \gMC \in\sMC_{\nmoment}^{\markk}} \prod_{\gC' }
    \sum_{\mathbf{x} \in D_{\gC'}} \alpha_{\gC'} f_{\gC'} \\ \nonumber & = \sum_{\pi}
    \prod\limits_{b \in \pi} \sum_{\gC' \in \sC_{b}^{\markk} }  \sum_{\mathbf{x} \in
    D_{\gC'}}  \alpha_{\gC'} f_{\gC'}                     \\ \nonumber & = \sum_{\pi} \prod\limits_{b
           \in \pi} \tilde{K}_{|b|},\end{align} where $\pi$ runs over all set-partitions
of the set $\{\underline{1}, \dots,  \underline{\nmoment}\}$, and $b$ runs over
all blocks of $\pi$. Above $\tilde{K}_j$ is defined as %
\begin{align} \label{eq:cumulant-posredni}
    \tilde{K}_j: & = \sum_{\gC \in
    \sC_j^{\markk} } \sum_{\mathbf{x} \in D_{\gC}} \alpha_{\gC} f_{\gC}                         \\ \nonumber
                 & = \frac{1}{j} \sum_{\gC \in \sC_j^{\markk} }   \sum_{\mathbf{x} \in D_{\gC}}
    \left( \prod\limits_{ v \in V_{\gC}} \mu_{\Lambda}(x_v) \right)
    (-1)^{|B_{\gC}|-1} f_{\gC}.\end{align} %

\medskip

In our context, the moment-cumulant formula \eqref{eq:rekurencja} takes the form:
\[ M_{\nmoment}=\sum_{\pi \in \Pi_{\nmoment}}\prod_{b\in \pi}K_{|b|}, \]
where $\pi$ ranges over set-partitions, and $b$ iterates through all blocks of
$\pi$. This formula can be interpreted as a system of algebraic equations for
the unknowns $(K_{\nmoment})$. The system exhibits an upper-triangular
structure, allowing us to express the $\nmoment$-th cumulant $K_{\nmoment}$ as
the sum of the moment $M_{\nmoment}$ and a complicated polynomial in $K_1, \dots,
    K_{\nmoment-1}$. This structure enables a recursive solution method,
guaranteeing a unique solution. Equation \eqref{eq:moment-cumulant2}
demonstrates that the sequence $(\tilde{K}_{\nmoment})$ satisfies this system of
equations. Given the uniqueness of the solution, we can conclude that the
cumulant
\[K_{\nmoment}=\tilde{K}_{\nmoment}\] is given by
\eqref{eq:cumulant-posredni} after the substitution $j=\nmoment$.

\subsection{Caterpillar graphs with labeled vertices}
\label{sec:cat-labeled}

\begin{figure}
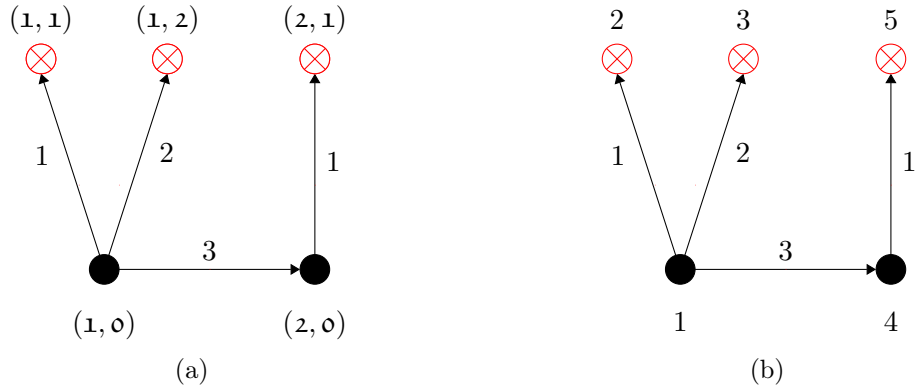

    \subfile{figures/FIGURE-tagged-and-labeled.tex}

    \caption{ \protect\subref{fig:labelb} A caterpillar graph with tagged vertices.
        \protect\subref{fig:labeld}~The same graph with labeled vertices. The labels
        belong to the set $\{ 1, 2, 3, 4, 5\}$.} \label{fig:labelowanie}
\end{figure}

We say that a connected, weighted, directed graph $\gC$ is a \emph{caterpillar
    graph with $k$ labeled vertices} if its vertex set is ${1, \dots, k}$, the edge
weights satisfy the convention \eqref{eq:labeling-compatible}, and there exists
a way to tag the vertices such that $\gC$ becomes an element of
$\sMC_{\nmoment}^{\tagg}$ (see \cref{sec:multi-caterpillar-formal}). An example
of a caterpillar graph with labeled vertices is shown in \cref{fig:labeld}. The
set of caterpillar graphs with $k$ labeled vertices is denoted by
$\sC^{\labeled}_{\nmoment}$. This definition may seem abstract, so we provide an
alternative description below.

\medskip

For any connected graph $\gMC \in
    \sMC_{\nmoment}^{\tagg}$, there is a unique way to label the vertices to satisfy
the requirement \eqref{eq:labeling-compatible}, as follows. Start by assigning
the number $1$ to the unique black vertex with no incoming edges. Then, in the
order given by the edge weights, number all endpoints of the edges outgoing from
this vertex with successive natural numbers. Repeat the process at the unique
black endpoint of an edge outgoing from vertex $1$, and continue until all black
vertices are visited. In this way, for any edge $e = (v_1, v_2)$, its weight is
equal to the difference of the labels of the endpoints: $\weight(e) = v_2 -
    v_1$. The result is clearly an element of $\sC^{\labeled}_{\nmoment}$, and each
element of $\sC^{\labeled}_{\nmoment}$ can be obtained in this manner.

\begin{figure}
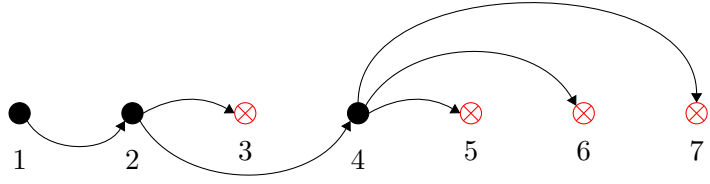

    \subfile{figures/FIGURE-red-black-tree} \caption{An example of a caterpillar
        graph with $\nmoment = 7$ labeled vertices. Using the notation from
        \eqref{eq:black}, we have $\ell = 3$ black vertices with $b_1 = 1$, $b_2 = 2$,
        and $b_3 = 4$. Additionally, we use the convention that $b_4 = 8$. The weights
        of the edges are not shown for clarity.}

    \label{fig:drzewocc}
\end{figure}

The above procedure shows that the elements of $\sC^{\labeled}_{\nmoment}$ can
be characterized as follows. For each $\gC \in \sC^{\labeled}_{\nmoment}$ with
the set of black vertices
\begin{equation}
    \label{eq:black}
    B_{\gC} = \{ b_1, \dots, b_\ell \} \subseteq \{ 1, \dots, \nmoment \}, \qquad b_1 < \dots < b_\ell
\end{equation}
we have $\ell \geq 1$ and $b_1 = 1$. We will use the convention that $b_{\ell+1}
    = \nmoment + 1$. The black vertices are connected by a series of directed edges:
\[ (b_1, b_2), \quad (b_2, b_3), \quad \dots, \quad (b_{\ell-1}, b_{\ell}).\]
Additionally, each black vertex $b_i$ (with $i \in \{1, \dots, \ell\}$) is
connected to the red vertices $b_i + 1, b_i + 2, \dots, b_{i+1} - 1$ that
immediately follow it by a collection of directed edges:
\[
    (b_i, b_i + 1), \quad (b_i, b_i + 2), \quad \dots, \quad (b_i, b_{i+1} - 1), \]
as illustrated in \cref{fig:drzewocc}.

In particular, since the structure of a caterpillar graph with labeled vertices
is determined by its set of black vertices, it follows that
$|\sC^{\labeled}_{\nmoment}| = 2^{\nmoment - 1}$.

\subsection{The second double counting. The second formula for the cumulants}
\label{sec:second-formula-cumulants}

We continue the discussion from \cref{sec:first-cumulant-formula} and revisit
the formula \eqref{eq:cumulant-posredni} for the cumulant $K_{\nmoment}$. As we
already mentioned, the connected graph $\gC \in \sC_\nmoment^{\markk}$ can be
labeled in a unique way so that it becomes an element of
$\sC_\nmoment^{\labeled}$. On the other hand, for each graph $\gC \in
    \sC^{\labeled}_{\nmoment}$, there exist $\nmoment !$ ways to mark the vertices
so that the outcome is a caterpillar graph with $\nmoment$ marked vertices. In
this way we proved the following intermediate result.

\begin{corollary}
    \label{cor:wzorek-na-kumulante}

    Let $\Lambda$ be an interlacing sequence with a generic set of concave corners.
    For each $u_0 \in \R$ the $n$-th formal cumulant considered in
    \cref{sec:cumulants-interlacing} is given by the following sum over caterpillar
    graphs with $\nmoment$~labeled vertices
    \begin{equation}
        \label{eq:kumulanta-caterpillar}
        K_{\nmoment} =  (\nmoment-1) !
        \sum_{\gC \in \sC^{\labeled}_{\nmoment} } \sum_{\mathbf{x} \in D_{\gC}}
        (-1)^{|\sB_{\gC}|-1}    f_{\gC} \prod\limits_{ j \in \{1,\dots,k\} }
        \mu_{\Lambda}(x_j).
    \end{equation}
\end{corollary}
For example, for $k=2$ we
obtain $$K_2=\sum_{x_1 \leq u_0}\sum_{x_2} \frac{
        \mu_{\Lambda}(x_1)\mu_{\Lambda}(x_2) }{x_2-x_1+1} - \sum_{x_1 \leq u_0}\sum_{x_2
        \leq u_0} \frac{ \mu_{\Lambda}(x_1)\mu_{\Lambda}(x_2) }{x_2-x_1+1}.$$ The first
summand corresponds to the caterpillar graph shown of
\cref{fig:lllterm}, and the second summand corresponds to the
caterpillar graph shown of \cref{fig:rrrterm}.

\begin{figure}
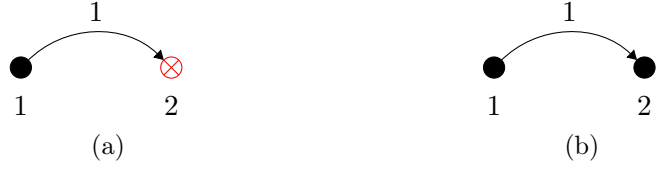

    \subfile{figures/FIGURE-caterrpillar-one-two.tex}
    \caption{
        \protect\subref{fig:lllterm} Caterpillar graph with one black and one red
        vertex.
        \protect\subref{fig:rrrterm}~Caterpillar graph with two black
        vertices.
    }
    \label{fig:pierwszedwa}
\end{figure}

\subsection{Sum over non-crossing alternating trees}
\label{sec:sum-ncat}

\begin{proposition} \label{cor:wzorek-na-kumulante2}
    Let $\Lambda$ be an
    interlacing sequence with a generic set of concave corners. For each $u_0 \in
        \R$ the $\nmoment$-th formal cumulant considered in
    \cref{sec:cumulants-interlacing} is given by the following sum over noncrossing
    alternating trees
    \begin{equation}
        \label{eq:kumulanta-trees}
        K_{\nmoment} =
        (\nmoment-1) ! \sum_{\gT \in \sT_{\nmoment} } \sum_{\mathbf{x} \in D_{\gT}}
        (-1)^{|\sB_{\gT}|-1}    f_{\gT} \prod\limits_{ j \in \{1,\dots,k\} }
        \mu_{\Lambda}(x_j).
    \end{equation}
\end{proposition}
\begin{proof}
    In
    \eqref{eq:kumulanta-caterpillar} we can reverse the order of the sums and write
    \[    K_{\nmoment} = -(\nmoment-1) ! \sum_{x_1,\dots,x_\nmoment \in \Deco}
        \mathfrak{C}_\nmoment(x_1,\dots,x_\nmoment) \prod\limits_{ j \in \{1,\dots,k\} }
        \mu_{\Lambda}(x_j),  \] where \[    \mathfrak{C}_\nmoment(x_1,\dots,x_\nmoment)
        := \sum_{\substack{\gC \in \sC^{\labeled}_{\nmoment} \\ (x_1,\dots,x_\nmoment)
                \in D_{\gC}} }  (-1)^{|\sB_{\gC}|}    f_{\gC}(x_1,\dots,x_\nmoment). \]
    Similarly, the right-hand side of \eqref{eq:kumulanta-trees} can be written as
    \[ -(\nmoment-1) ! \sum_{x_1,\dots,x_\nmoment \in \Deco}
        \mathfrak{T}_\nmoment(x_1,\dots,x_\nmoment) \prod\limits_{ j \in \{1,\dots,k\} }
        \mu_{\Lambda}(x_j),  \] where \[    \mathfrak{T}_\nmoment(x_1,\dots,x_\nmoment)
        := \sum_{\substack{\gT \in \sT_{\nmoment} \\ (x_1,\dots,x_\nmoment) \in D_{\gT}}
        }  (-1)^{|\sB_{\gT}|}    f_{\gT}(x_1,\dots,x_\nmoment). \] As a side remark note
    that $\mathfrak{T}_\nmoment(x_1,\dots,x_\nmoment)$ is a quantity which (up to a
    scaling factor) is closely related to the random variable $Z$ from
    \cref{rem:cumulant-as-expected-value}. The result is a consequence of
    \cref{cor:katerpilary-na-drzewa} below.
\end{proof}

\begin{lemma} \label{cor:katerpilary-na-drzewa}
    With the above notations,
    \[
        \mathfrak{C}_\nmoment(x_1,\dots,x_\nmoment)
        =\mathfrak{T}_\nmoment(x_1,\dots,x_\nmoment) \]
    holds true for any $\nmoment\geq
        1$ and any $x_1,\dots,x_{\nmoment}\in\R$ for which the left-hand side of the
    equality does not involve division by zero.
\end{lemma}
\begin{proof}
    In the
    special case $k=1$ we have that the set of graphs $\sC^{\labeled}_1=\sT_1$ which
    contributes to $\mathfrak{C}_1(x_1)$, respectively to $\mathfrak{T}_1(x_1)$,
    consists of a single element depicted in \cref{fig:llterm}. Thus \[
        \mathfrak{C}_1(x_1)=\mathfrak{T}_1(x_1)= \begin{cases} -1 & \text{if } x_1 \leq
              u_0,                     \\ 0           & \text{if } x_1 > u_0.\end{cases} \]

    \medskip

    Let $\nmoment \geq 2$. In the case when $x_{k}\leq u_0$ we obtain that \[
        \mathfrak{T}_{x_1,\dots,x_k} = 0 \] because the rightmost vertex of any
    non-crossing alternating tree $\gT\in \sT_k$ is white thus
    $(x_1,\dots,x_\nmoment)$ is not a decoration of $\gT$ and the sum runs over the
    empty set.

    \medskip

    Let $\sC^{\labeled}(V)$ and $\sT(V)$ denote, respectively, the set of
    caterpillar graphs and the set of non-crossing alternating trees with vertex set
    $V$. In particular, $\sC^{\labeled}(\{1,\ldots,\nmoment\}) =
        \sC^{\labeled}_{\nmoment}$ and $\sT(\{1,\ldots,\nmoment\}) = \sT_{\nmoment}$ for
    any natural number $\nmoment$.

    Let $\gT \in \sT_\nmoment$ be a non-crossing alternating tree with $\nmoment$
    vertices. Clearly, $\gT$ contains the edge $e=(1,\nmoment)$ connecting the
    leftmost and rightmost vertices. Removing $e$ from $\gT$ results in two
    connected components: $\gT_1 \in \sT_{i-1}(\{1,\ldots,i-1\})$ and $\gT_2 \in
        \sT_{\nmoment-i}(\{i,\ldots,\nmoment\})$ for some $i \in \{2,\ldots,\nmoment\}$.
    If $\gT_2$ consists of a single (white) vertex, we change its color to black.
    Thus, each non-crossing alternating tree with vertex set $\{1,\ldots,\nmoment\}$
    decomposes uniquely into the edge $e$ and two non-crossing alternating trees
    $\gT_1$ and $\gT_2$ with vertex sets $\{1,\ldots,i-1\}$ and
    $\{i,\ldots,\nmoment\}$, respectively.

    For $i \neq \nmoment$, the tuple $(x_1,\ldots,x_\nmoment)$ is a decoration of $\gT$ if and
    only if $(x_1,\ldots,x_{i-1})$ is a decoration of $\gT_1$ and
    $(x_i,\ldots,x_{\nmoment})$ is a decoration of $\gT_2$. The case where
    $i=\nmoment$ and $\gT_2$ consists of a single vertex requires separate
    consideration. This leads to the following recurrence relation:
    \[
        \mathfrak{T}_{\nmoment}(x_1,\ldots,x_{\nmoment}) =
        \begin{cases}
            0                                                                                        & \text{if } x_{\nmoment} \leq u_0, \\[2ex]
            \sum\limits_{i=2}^{\nmoment-1} \frac{\mathfrak{T}_{i-1}(x_1,\ldots,x_{i-1}) \mathfrak{T}_{\nmoment-i+1}(x_i,\ldots,x_{\nmoment})}{x_{\nmoment}-x_1+\nmoment-1} +
            \frac{\mathfrak{T}_{\nmoment-1}(x_1,\ldots,x_{\nmoment-1})}{x_{\nmoment}-x_1+\nmoment-1} & \text{if } x_{\nmoment} > u_0.
        \end{cases}
    \]
    We will prove that the sequence of functions $\mathfrak{C}_{\nmoment}$ satisfies
    the same recurrence relation.

    \bigskip

    If $\nmoment\geq 2$ and $x_{\nmoment} \leq u_0$ then $
        \mathfrak{C}(x_1,\dots,x_\nmoment)=0$ %
    because we can pair caterpillar graphs from $\sC^{\labeled}_\nmoment$ into pairs
    that differ only in the color of the far-right vertex, and the contribution of
    each pair to the sum is zero.

    \medskip

    Let $k \geq 2$ be a natural number, let $x_{\nmoment} > u_0$ and let $\gC \in
        \sC^{\labeled}_{\nmoment}$ be a caterpillar graph. There is a unique path $(i_1,
        \dots, i_t)$ with $t\geq 2$ from the vertex $1$ to the vertex $\nmoment$ which
    means that $i_1=1$ and $i_t=\nmoment$, and $(i_1,i_2), (i_2,i_3), \dots,
        (i_{t-1},i_t) \in \sE_{\gC}$. %
    In the special case when $e=(1, \nmoment) \in E_{\gC}$, we have $t=2$ and
    $e_1=e$. Using the telescopic sum
    $$x_{\nmoment}-x_{1}+\nmoment-1=\sum\limits_{j=1}^{t-1}
        \left(x_{i_{j+1}}-x_{i_j} + i_{j+1}-i_j\right)$$ we obtain
    \begin{align*}
        f_{\gC}(x_1,\dots,x_\nmoment) & = \sum\limits_{j=1}^{t-1}\frac{
        x_{i_{j+1}}-x_{i_j} + i_{j+1}-i_j }{x_{\nmoment}-x_1+\nmoment-1}
        f_{\gC}(x_1,\dots,x_\nmoment)                                   \\      & =\sum\limits_{j=1}^{t-1}\frac{ f_{\gC \setminus
                           e_j}(x_1,\dots,x_\nmoment) }{x_{\nmoment}-x_1+\nmoment-1},
    \end{align*}
    where $
        \gC \setminus e_j $ denotes the graph $\gC$ with the edge $e_j$ removed.
    Therefore, $$ \mathfrak{C}(x_1,\dots,x_\nmoment) =\sum\limits_{\substack{ \gC
                \in \sC^{\labeled}_{\nmoment} \\ (x_1,\dots,x_\nmoment)\in D_{\gC} }}
        (-1)^{|B_{\gC}|}\sum\limits_{j=1}^{t-1}\frac{ f_{\gC \setminus e_j
                }(x_1,\dots,x_\nmoment) }{x_{\nmoment}-x_1+\nmoment-1}. $$

    In addition, every graph $ \gC \in \sC^{\labeled}_{\nmoment}$ after removing any
    edge $e_j$ splits in a unique way into the sum of two caterpillar graphs
    ${\gC}_1 \in \sC^{\labeled}_{i-1}\left(\{ 1, \dots, i-1 \} \right)$ and ${\gC}_2
        \in \sC^{\labeled}_{\nmoment-i+1}\left( \{ i, \dots, \nmoment \} \right) $ for
    some $i\in\{2,\dots,k\}$. %
    In the special case when $i=\nmoment$ and the graph $\gC_2$ consists a single
    (red) vertex, we change its color to black; this case will require separate
    analysis. In this way we can write $\mathfrak{C}(x_1,\dots,x_\nmoment)$ as a
    triple sum over all possible choices of the number $i$, over all graphs $\gC_1
        \in \sC^{\labeled}_{i-1}\left(\{ 1, \dots, i-1 \} \right)$ and over all graphs
    $\gC_2 \in \sC^{\labeled}_{\nmoment-i+1}\left( \{ i, \dots, \nmoment \}
        \right)$, i.e.,
    \begin{multline*}
        \mathfrak{C}(x_1,\dots,x_\nmoment)  = \\
        \sum_{i=2}^{\nmoment-1} \; \sum_{ \substack{ \gC_1 \in
        \sC^{\labeled}_{i-1}\left(\{ 1, \dots, i-1 \} \right) \\ (x_1,\dots,x_{i-1}) \in
        D_{\gC_1}}} \; \sum_{ \substack{ \gC_2 \in \sC^{\labeled}_{\nmoment-i+1}\left(
        \{ i, \dots, \nmoment \} \right) \\ (x_i,\dots,x_{\nmoment}) \in D_{\gC_2} }}
        \frac{(-1)^{|B_{\gC_1}|+|B_{\gC_2}|} f_{\gC_1} f_{\gC_2}
        }{x_{\nmoment}-x_1+\nmoment-1} \\ %
        \shoveright{+ \sum\limits_{ \substack{ \gC_1 \in \sC^{\labeled}_{i-1}\left(\{ 1,
                \dots, \nmoment-1 \} \right) \\ B_{\gC_1} \subseteq B }} \;
        \frac{(-1)^{|B_{\gC_1}|} f_{\gC_1}  }{x_{\nmoment}-x_1+\nmoment-1}    } %
        \\ = \sum\limits_{i=2}^{\nmoment-1} \frac{\mathfrak{C}_{i-1}(x_1, \dots, x_{i-1}
            ) \mathfrak{C}_{\nmoment-i+1}(x_i, \dots, x_{\nmoment} )
        }{x_{\nmoment}-x_1+\nmoment-1} +  \frac{\mathfrak{C}_{\nmoment-1}(x_1, \dots,
            x_{\nmoment-1} ) }{x_{\nmoment}-x_1+\nmoment-1}
    \end{multline*}
    if $x_\nmoment>
        u_0$.

    The sequences of rational functions $\mathfrak{C}_\nmoment$ and
    $\mathfrak{T}_\nmoment$ satisfy the same recurrence relation and have the same
    initial condition, which completes the proof.
\end{proof}

\subsection{Proof of \cref{thm:wzormikolaja}}
\label{sec:proof-limit}

\begin{proof}[Proof of \cref{thm:wzormikolaja}]
    When the number $u_0$ is not an integer, we apply \cref{cor:regularize2} and
    evaluate the cumulant $K_{\nmoment} (\Lambda^\epsilon, u_0)$ using
    \cref{cor:wzorek-na-kumulante2}.

    However, when the number $u_0$ is an integer, the theorem is satisfied for any
    number $u \in (u_0, u_0+1)$, as shown above. Since $F_T(u_0)$ is a
    right-continuous function, then \[F_T(u_0)=\lim\limits_{u \to u_0} F_T(u),\] and
    \cref{thm:wzormikolaja} also holds for the number $u_0$.
\end{proof}

%% file: figures/FIGURE-Pieri.tex
	        \begin{tikzpicture}[scale=0.7,rotate=45]
		\coordinate (p0) at (6,0);
		\coordinate (p1) at (6,2);
		\coordinate (p2) at (5,2);
		\coordinate (p3) at (5,4);
		\coordinate (p4) at (3,4);
		\coordinate (p5) at (3,6);
		\coordinate (p6) at (2,6);
		\coordinate (p7) at (2,8);
		\coordinate (p8) at (0,8);
		\coordinate (p9) at (0,12);
		\coordinate (b1) at (6.5,0.5);
		\coordinate (b2) at (5.5,2.5);
		\coordinate (b3) at (3.5,4.5);
		\coordinate (b4) at (2.5,6.5);
		\coordinate (b5) at (.5,8.5);
		\coordinate (r2) at (5.5,3.5);
		\coordinate (r3) at (3.5,5.5);
		\coordinate (r4) at (0.5,9.5);
		\coordinate (r5) at (0.5,10.5);
		\coordinate (tam) at (0.65,-0.35);
		\coordinate (tam2) at (0.4,0.4);
		\coordinate (gora) at (0,1.2);
		\coordinate (dol) at (0,-0.9);
		\coordinate (lewo) at (-.5,0.1);
		\coordinate (prawo) at (.5,0.1);
		\coordinate (lewodol) at (-.5,-0.7);
		\coordinate (prawodol) at (.5,-0.7);
		
		\draw[draw=orange,fill=orange!10] (p0) -- (p1) -- (p2) -- (p3) -- (p4) -- (p5) -- (p6) -- (p7) -- (p8) -- (p9) -- (0,0) -- cycle;

		\fill [fill=blue!20] (p4) rectangle +(1,1);
		\fill [fill=blue!20] (p8) rectangle +(1,3);
		\fill [fill=blue!20] (p0) rectangle +(1,2);
        
		\draw [blue] (p4) grid +(1,1);
        \draw [blue] (p8) grid +(1,3);
        \draw [blue] (p0) grid +(1,2);

		\draw[very thick, orange] (8,0) -- (p0) -- (p1) -- (p2) -- (p3) -- (p4) -- (p5) -- (p6) -- (p7) -- (p8) -- (p9) ;
		\draw[->] (-5.8,5.8) -- (3.8,-3.8) node[anchor=west] {$u$};
		\draw[dotted,blue,very thick] (p8) -- (-8/2,8/2) node[anchor=north]{\textcolor{black}{$x_3$}};
		\draw[dotted,blue,very thick] (p4) -- (-1/2,1/2) node[anchor=north]{\textcolor{black}{$x_2$}};
		\draw[dotted,blue,very thick] (p0) -- (6/2,-6/2) node[anchor=north]{\textcolor{black}{$x_1$}};
	\end{tikzpicture}
	

%% file: figures/FIGURE-spine.tex
		\begin{center}
		\begin{tikzpicture}[scale=1.0]
			\punkty	
			\wierzcholek{$1.5*(x1)$}{black}{3}{dol}{black}
			\wierzcholek{$1.5*(x2)$}{black}{5}{dol}{black}
			\wierzcholek{$1.5*(x3)$}{black}{1}{dol}{black}
			\wierzcholek{$1.5*(x4)$}{black}{4}{dol}{black}
			\wierzcholek{$1.5*(x5)$}{black}{6}{dol}{black}
			\wierzcholek{$1.5*(x6)$}{black}{2}{dol}{black}
			\strzalka{$1.5*(x1)$}{$1.5*(x2)+(-3pt,0)$}{black}
			\strzalka{$1.5*(x2)$}{$1.5*(x3)+(-3pt,0)$}{black}
			\strzalka{$1.5*(x3)$}{$1.5*(x4)+(-3pt,0)$}{black}
			\strzalka{$1.5*(x4)$}{$1.5*(x5)+(-3pt,0)$}{black}
			\strzalka{$1.5*(x5)$}{$1.5*(x6)+(-3pt,0)$}{black}
		\end{tikzpicture}
	\end{center}
	

%% file: figures/FIGURE-multi-spine.tex
		\begin{center}
		\begin{tikzpicture}[scale=1.0]
			\punkty	
			\wierzcholek{$1.5*(x1)$}{black}{4}{dol}{black}
			\wierzcholek{$1.5*(x2)$}{black}{7}{dol}{black}
			\wierzcholek{$1.5*(x3)$}{black}{1}{dol}{black}
			\wierzcholek{$1.5*(x4)$}{black}{6}{dol}{black}
			\wierzcholek{$1.5*(x5)$}{black}{2}{dol}{black}
			\wierzcholek{$1.5*(x6)$}{black}{8}{dol}{black}
			\wierzcholek{$1.5*(x7)$}{black}{5}{dol}{black}
			\wierzcholek{$1.5*(x8)$}{black}{9}{dol}{black}
			\wierzcholek{$1.5*(x9)$}{black}{3}{dol}{black}
			\strzalka{$1.5*(x1)$}{$1.5*(x2)+(-3pt,0)$}{black}
			\strzalka{$1.5*(x2)$}{$1.5*(x3)+(-3pt,0)$}{black}
			\strzalka{$1.5*(x5)$}{$1.5*(x6)+(-3pt,0)$}{black}
			\strzalka{$1.5*(x7)$}{$1.5*(x8)+(-3pt,0)$}{black}
			\strzalka{$1.5*(x8)$}{$1.5*(x9)+(-3pt,0)$}{black}
		\end{tikzpicture}
	\end{center}
	

%% file: figures/FIGURE-graph-on-diagram.tex
	        \begin{tikzpicture}[scale=0.9,rotate=45]
		\coordinate (p0) at (6,0);
		\coordinate (p1) at (6,2);
		\coordinate (p2) at (5,2);
		\coordinate (p3) at (5,4);
		\coordinate (p4) at (3,4);
		\coordinate (p5) at (3,6);
		\coordinate (p6) at (2,6);
		\coordinate (p7) at (2,8);
		\coordinate (p8) at (0,8);
		\coordinate (p9) at (0,12);
		\coordinate (b1) at (6.5,0.5);
		\coordinate (b2) at (5.5,2.5);
		\coordinate (b3) at (3.5,4.5);
		\coordinate (b4) at (2.5,6.5);
		\coordinate (b5) at (.5,8.5);
		\coordinate (r2) at (5.5,3.5);
		\coordinate (r3) at (3.5,5.5);
		\coordinate (r4) at (0.5,9.5);
		\coordinate (r5) at (0.5,10.5);
		\coordinate (tam) at (0.65,-0.35);
		\coordinate (tam2) at (0.4,0.4);
		\coordinate (gora) at (0,1.2);
		\coordinate (dol) at (0,-0.9);
		\coordinate (lewo) at (-.5,0.1);
		\coordinate (prawo) at (.5,0.1);
		\coordinate (lewodol) at (-.5,-0.7);
		\coordinate (prawodol) at (.5,-0.7);
		
		\draw[orange,fill=orange!15]
			(p0) -- (p1) -- (p2) -- (p3) -- (p4) -- (p5) -- (p6) -- (p7) -- (p8) -- (p9) -- (0,0) -- cycle;
		\draw [blue] (p6) grid +(1,1);
		\draw [blue] (p4) grid +(1,2);
		\draw [blue] (p2) grid +(1,2);
		\draw [blue] (p8) grid +(1,3);
		\draw [blue] (p0) grid +(1,1);
		\draw[very thick, orange] (8,0) -- (p0) -- (p1) -- (p2) -- (p3) -- (p4) -- (p5) -- (p6) -- (p7) -- (p8) -- (p9) ;
		\draw[->] (-5.8,5.8) -- (3.8,-3.8) node[anchor=west] {$u$};
		\draw[dotted,blue,very thick] (p8) -- (-8/2,8/2);
		\draw[dotted,blue,very thick] (p6) -- (-4/2,4/2);
		\draw[dotted,blue,very thick] (p4) -- (-1/2,1/2);
		\draw[dotted,blue,very thick] (p2) -- (3/2,-3/2);
		\draw[dotted,blue,very thick] (p0) -- (6/2,-6/2);
		\foreach \u in {-11,...,7}
		{\draw (\u/2,-\u/2)  +(-1pt,-1pt) node[anchor=north] {\tiny $\u$} -- +(1pt,1pt); }
		\lstrzalka{$(b1)+(-4pt, 0pt)$}{$(b5)+(0pt, -4 pt)$}{45}{black}
		\strzalka{b2}{$(r2)+(0 , -4 pt)$}{black}
		\strzalka{b3}{$(r3)+(0,-4pt)$}{black}
		\strzalka{b5}{$(r4)+(0,-4pt)$}{black}
		\pstrzalka{b5}{$(r5)+(2.83 pt,-2.83 pt)$}{30}{black}
		\strzalka{b4}{$(b3)+(-2.83 pt, 2.83 pt)$}{black}
		\strzalka{b5}{$(b4)+(-2.83 pt,2.83 pt)$}{black}
		\bwierzcholek{b1}{black}{\mathfrak{(1,0)}}{tam}{black}
		\bwierzcholek{b2}{black}{\mathfrak{(2,0)}}{tam}{black}
		\bwierzcholek{b3}{black}{\mathfrak{(3,0)}}{tam}{black}
		\bwierzcholek{b4}{black}{\mathfrak{(4,0)}}{tam2}{black}
		\bwierzcholek{b5}{black}{\mathfrak{(5,0)}}{tam}{black}
		\wierzcholekcbok{r2}{red}{\mathfrak{(2,1)}}{tam}{black}
		\wierzcholekcbok{r3}{red}{\mathfrak{(3,1)}}{tam}{black}
		\wierzcholekcbok{r4}{red}{\mathfrak{(5,1)}}{tam}{black}
		\wierzcholekcbok{r5}{red}{\mathfrak{(5,2)}}{tam}{black}
	\end{tikzpicture}
	

%% file: figures/FIGURE-graph-without-diagram.tex
	        \begin{tikzpicture}[scale=1.3]
		\punkty
		\strzalka{x9}{$(y9)+(0, -4pt)$}{black}
		\strzalka{x7}{$(y7)+(0, -4pt)$}{black}
		\strzalka{x3}{$(yy2)+(0, -4pt)$}{black}
		\strzalka{x1}{$(x3)+(-4pt, 0)$}{black}
		\strzalka{x3}{$(x5)+(-4pt, 0)$}{black}	
		\strzalka{x5}{$(x7)+(-4pt, 0)$}{black}		
		\strzalka{x3}{$(yy3)+(0, -4pt)$}{black}		
		\wierzcholekcpusty{x1}{red}{\mathfrak{(1,0)}}{dol}{black}
		\wierzcholekcpusty{x3}{red}{\mathfrak{(5,0)}}{dol}{black}	
		\wierzcholekcpusty{x5}{red}{\mathfrak{(4,0)}}{dol}{black}	
		\wierzcholekcpusty{x7}{red}{\mathfrak{(3,0)}}{dol}{black}	
		\wierzcholekcpusty{x9}{red}{\mathfrak{(2,0)}}{dol}{black}		
		\wierzcholek{x1}{black}{}{lewo}{blue}
		\wierzcholek{x3}{black}{}{lewo}{blue}
		\wierzcholek{x5}{black}{}{lewo}{blue}
		\wierzcholek{x7}{black}{}{lewo}{blue}
		\wierzcholek{x9}{black}{}{lewo}{blue}
		\wierzcholekpusty{y1}{white}{}{prawo}{blue}
		\wierzcholekpusty{y3}{white}{}{prawo}{blue}
		\wierzcholekc{y9}{red}{\mathfrak{(2,1)}}{$3.0*(gora)$}{black}
		\wierzcholekc{y7}{red}{\mathfrak{(3,1)}}{$3.0*(gora)$}{black}
		\wierzcholekc{yy2}{red}{\mathfrak{(5,1)}}{$3.0*(gora)$}{black}
		\wierzcholekc{yy3}{red}{\mathfrak{(5,2)}}{$3.0*(gora)$}{black}
		\wierzcholekpusty{$.6*(x3)+.4*(yy2)$}{black}{1}{$(lewo)$}{black}
		\wierzcholekpusty{$.6*(x3)+.4*(yy3)$}{black}{2}{$(prawo)$}{black}
		\wierzcholekpusty{$.6*(x9)+.4*(y9)$}{black}{1}{$.6*(lewo)$}{black}
		\wierzcholekpusty{$.6*(x7)+.4*(y7)$}{black}{1}{$.6*(lewo)$}{black}
		\wierzcholekpusty{x2}{red}{1}{gora}{black}
		\wierzcholekpusty{x4}{red}{3}{gora}{black}
		\wierzcholekpusty{x6}{red}{1}{gora}{black}
	\end{tikzpicture}
	

%% file: figures/FIGURE-tagged-and-marked.tex
	    \subfloat[]
	{
		\begin{tikzpicture}[scale=1.4]
			\punkty
			\strzalka{x1}{$(x3)+(-4pt,0)$}{black}
			\strzalka{x3}{$(y3)+(0,-4pt)$}{black}
			\strzalka{x1}{$(y1)+(0,-4pt)$}{black}		
			\wierzcholekcpusty{x1}{red}{\mathfrak{(1,0)}}{$1.05*(dol)$}{black}
			\wierzcholekcpusty{x3}{red}{\mathfrak{(2,0)}}{$1.05*(dol)$}{black}		
			\wierzcholek{x1}{black}{}{$.7*(lewo)+.7*(lewodol)+(0, .2)$}{blue}
			\wierzcholek{x3}{black}{}{$.7*(prawo)+.7*(prawodol)+(0, .2)$}{blue}
			\wierzcholekpusty{y1}{white}{}{$.7*(lewo)+.7*(lewodol)+(0, .2)$}{blue}
			\wierzcholekpusty{y3}{white}{}{$.7*(prawo)+.7*(prawodol)+(0, .2)$}{blue}
			\wierzcholekc{y1}{red}{\mathfrak{(1,1)}}{$3.0*(gora)$}{black}
			\wierzcholekc{y3}{red}{\mathfrak{(2,1)}}{$3.0*(gora)$}{black}
			
			\wierzcholekpusty{$.6*(x3)+.4*(y3)$}{red}{1}{$.5*(lewo)$}{black}
			\wierzcholekpusty{x2}{red}{2}{$.1*(gora)$}{black}
			\wierzcholekpusty{$.6*(x1)+.4*(y1)$}{red}{1}{$.5*(prawo)$}{black}
		\end{tikzpicture}
		\label{fig:przyka}
	}
	\hfill
	\subfloat[]
	{
		\begin{tikzpicture}[scale=1.4]
			\punkty
			\strzalka{x1}{$(x3)+(-4pt,0)$}{black}
			\strzalka{x3}{$(y3)+(0,-4pt)$}{black}
			\strzalka{x1}{$(y1)+(0,-4pt)$}{black}		
			\wierzcholekcpusty{x1}{red}{}{dol}{black}
			\wierzcholekcpusty{x3}{red}{}{dol}{black}		
			\wierzcholek{x1}{black}{\underline{4}}{$.7*(lewo)+.7*(lewodol)+(0, .2)$}{blue}
			\wierzcholek{x3}{black}{\underline{1}}{$.7*(prawo)+.7*(prawodol)+(0, .2)$}{blue}
			\wierzcholekpusty{y1}{white}{\underline{2}}{$.55*(lewo)+.85*(lewodol)+(0, .2)$}{blue}
			\wierzcholekpusty{y3}{white}{\underline{3}}{$.55*(prawo)+.85*(prawodol)+(0, .2)$}{blue}
			\wierzcholekc{y1}{red}{}{gora}{black}
			\wierzcholekc{y3}{red}{}{gora}{black}
			
			\wierzcholekpusty{$.6*(x3)+.4*(y3)$}{red}{1}{$.5*(lewo)$}{black}
			\wierzcholekpusty{x2}{red}{2}{$.1*(gora)$}{black}
			\wierzcholekpusty{$.6*(x1)+.4*(y1)$}{red}{1}{$.5*(prawo)$}{black}
		\end{tikzpicture}
		\label{fig:przykb}
	}
	\hfill
	\subfloat[]
	{
		\begin{tikzpicture}[scale=1.4]
			\punkty
			\strzalka{x1}{$(x3)+(-4pt,0)$}{black}
			\strzalka{x3}{$(y3)+(0,-4pt)$}{black}
			\strzalka{x1}{$(y1)+(0,-4pt)$}{black}		
			\wierzcholekcpusty{x1}{red}{\mathfrak{(1,0)}}{$1.05*(dol)$}{black}
			\wierzcholekcpusty{x3}{red}{\mathfrak{(2,0)}}{$1.05*(dol)$}{black}
			\wierzcholek{x1}{black}{\underline{4}}{$.7*(lewo)+.7*(lewodol)+(0, .2)$}{blue}
			\wierzcholek{x3}{black}{\underline{1}}{$.7*(prawo)+.7*(prawodol)+(0, .2)$}{blue}
			\wierzcholekpusty{y1}{white}{\underline{2}}{$.55*(lewo)+.85*(lewodol)+(0, .2)$}{blue}
			\wierzcholekpusty{y3}{white}{\underline{3}}{$.55*(prawo)+.85*(prawodol)+(0, .2)$}{blue}
			\wierzcholekc{y1}{red}{\mathfrak{(1,1)}}{$3.0*(gora)$}{black}
			\wierzcholekc{y3}{red}{\mathfrak{(2,1)}}{$3.0*(gora)$}{black}
			
			\wierzcholekpusty{$.6*(x3)+.4*(y3)$}{red}{1}{$.5*(lewo)$}{black}
			\wierzcholekpusty{x2}{red}{2}{$.1*(gora)$}{black}
			\wierzcholekpusty{$.6*(x1)+.4*(y1)$}{red}{1}{$.5*(prawo)$}{black}
		\end{tikzpicture}
		\label{fig:przykc}
	}
	

%% file: figures/FIGURE-tagged-and-labeled.tex
	\subfloat[]
	{
		\begin{tikzpicture}[scale=1.4]
			\punkty
			\strzalka{x1}{$(x3)+(-4pt,0)$}{black}
			\strzalka{x3}{$(y3)+(0,-4pt)$}{black}
			\strzalka{x1}{$(y00)+(0,-4pt)$}{black}
			\strzalka{x1}{$(y11)+(0,-4pt)$}{black}				

			\wierzcholekcpusty{x1}{red}{\mathfrak{(1,0)}}{dol}{black}
			\wierzcholekcpusty{x3}{red}{\mathfrak{(2,0)}}{dol}{black}
			\wierzcholekc{y00}{red}{\mathfrak{(1,1)}}{$3.0*(gora)$}{black}
			\wierzcholekc{y11}{red}{\mathfrak{(1,2)}}{$3.0*(gora)$}{black}
			\wierzcholekc{y3}{red}{\mathfrak{(2,1)}}{$3.0*(gora)$}{black}
			
			\wierzcholekcpusty{x1}{red}{}{dol}{black}
			\wierzcholekcpusty{x3}{red}{}{dol}{black}		
			\wierzcholek{x1}{black}{}{$.7*(lewo)+.7*(lewodol)+(0, .2)$}{blue}
			\wierzcholek{x3}{black}{}{$.7*(prawo)+.7*(prawodol)+(0, .2)$}{blue}
			\wierzcholekpusty{y1}{white}{}{$.7*(lewo)+.7*(lewodol)+(0, .2)$}{blue}
			\wierzcholekpusty{y3}{white}{}{$.7*(prawo)+.7*(prawodol)+(0, .2)$}{blue}
			\wierzcholekc{y00}{red}{}{gora}{black}
			\wierzcholekc{y11}{red}{}{gora}{black}
			\wierzcholekc{y3}{red}{}{gora}{black}
			
			\wierzcholekpusty{$.6*(x3)+.4*(y3)$}{red}{1}{$.5*(prawo)$}{black}
			\wierzcholekpusty{x2}{red}{3}{$.1*(gora)$}{black}
			\wierzcholekpusty{$.5*(x1)+.4*(y11)$}{red}{2}{$1.3*(prawo)$}{black}
			\wierzcholekpusty{$.6*(x1)+.4*(y00)$}{red}{1}{$(lewo)$}{black}
		\end{tikzpicture}
		\label{fig:labelb}
	}
	\hfill
	\subfloat[]
	{
		\begin{tikzpicture}[scale=1.4]
			\punkty
			\strzalka{x1}{$(x3)+(-4pt,0)$}{black}
			\strzalka{x3}{$(y3)+(0,-4pt)$}{black}
			\strzalka{x1}{$(y00)+(0,-4pt)$}{black}
			\strzalka{x1}{$(y11)+(0,-4pt)$}{black}			
			\wierzcholekcpusty{x1}{red}{\mathfrak{}}{dol}{black}
			\wierzcholekcpusty{x3}{red}{\mathfrak{}}{dol}{black}		
			\wierzcholek{x1}{black}{1}{$.9*(dol)$}{black}
			\wierzcholek{x3}{black}{4}{$.9*(dol)$}{black}
			\wierzcholekpusty{y11}{white}{3}{$4.0*(gora)$}{black}
			\wierzcholekpusty{y00}{white}{2}{$4.0*(gora)$}{black}
			\wierzcholekpusty{y3}{white}{5}{$4.0*(gora)$}{black}
			\wierzcholekc{y00}{red}{\mathfrak{}}{gora}{black}
			\wierzcholekc{y11}{red}{}{gora}{black}
			\wierzcholekc{y3}{red}{}{gora}{black}

			\wierzcholekpusty{$.6*(x3)+.4*(y3)$}{red}{1}{$.5*(prawo)$}{black}
			\wierzcholekpusty{x2}{red}{3}{$.1*(gora)$}{black}
			\wierzcholekpusty{$.5*(x1)+.4*(y11)$}{red}{2}{$1.3*(prawo)$}{black}
			\wierzcholekpusty{$.6*(x1)+.4*(y00)$}{red}{1}{$(lewo)$}{black}
		\end{tikzpicture}
		\label{fig:labeld}
	}
	

%% file: figures/FIGURE-red-black-tree.tex
		\begin{center}
		\begin{tikzpicture}[scale=1.0]
			\punkty	
			\wierzcholek{$1.5*(x1)$}{black}{1}{$1.15*(dol)$}{black}
			\wierzcholek{$1.5*(x2)$}{black}{2}{$1.15*(dol)$}{black}
			\wierzcholekc{$1.5*(x3)$}{black}{3}{dol}{black}
			\wierzcholek{$1.5*(x4)$}{black}{4}{$1.15*(dol)$}{black}
			\wierzcholekc{$1.5*(x5)$}{black}{5}{dol}{black}
			\wierzcholekc{$1.5*(x6)$}{black}{6}{dol}{black}
			\wierzcholekc{$1.5*(x7)$}{black}{7}{dol}{black}
			\pstrzalka{$1.5*(x1)+(2.83 pt, -2.83 pt)$}{$1.5*(x2)+(-2.83 pt, -2.83 pt )$}{60}{black}
			\pstrzalka{$1.5*(x2)+(2.83 pt, -2.83 pt)$}{$1.5*(x4)+(-2.83 pt, -2.83 pt )$}{60}{black}
			\lstrzalka{$1.5*(x2)+(4 pt, 0)$}{$1.5*(x3)+(-4 pt, 0 pt)$}{30}{black}
			\lstrzalka{$1.5*(x4)+(4 pt, 0)$}{$1.5*(x5)+(-4 pt, 0 pt)$}{30}{black}
			\lstrzalka{$1.5*(x4)+(2.83 pt, 2.83 pt )$}{$1.5*(x6)+(-2.83 pt, 2.83 pt)$}{60}{black}
			\lstrzalka{$1.5*(x4)+(0, 4 pt)$}{$1.5*(x7)+(0 pt, 4 pt)$}{90}{black}
		\end{tikzpicture}
	\end{center}

%% file: figures/FIGURE-caterrpillar-one-two.tex
			\begin{center}
			\subfloat[]
			{
				\begin{tikzpicture}[scale=1.0]
					\punkty	
					\lstrzalka{$(x1)+(2.12pt,2.12pt)$}{$(x3)+(-2.83pt,2.83pt)$}{45}{black}
					\wierzcholek{x1}{black}{1}{dol}{black}
					\wierzcholekpusty{x2}{white}{1}{$10*(gora)$}{black}
					\wierzcholekc{x3}{white}{2}{dol}{black}
				\end{tikzpicture}
				\label{fig:lllterm}
			}
			\qquad
			\qquad
			\qquad
			\qquad
			\qquad
			\subfloat[]
			{
				\begin{tikzpicture}[scale=1.0]
					\punkty	
					\lstrzalka{$(x1)+(2.12pt,2.12pt)$}{$(x3)+(-2.83pt,2.83pt)$}{45}{black}
					\wierzcholek{x1}{black}{1}{dol}{black}
					\wierzcholekpusty{x2}{white}{1}{$10*(gora)$}{black}
					\wierzcholek{x3}{black}{2}{dol}{black}
				\end{tikzpicture}
				\label{fig:rrrterm}
			}
		\end{center}